\newcount\Comments  
\documentclass[letterpaper,10pt]{scrartcl}

\usepackage[utf8]{inputenc}
\usepackage[english]{babel}
\usepackage[automark]{scrlayer-scrpage}
\usepackage{amsmath,amstext,amsbsy,amsopn,amscd,amsxtra,upref,amssymb}
\usepackage{amsfonts,bm}
\usepackage{amsthm}
\usepackage{color}
\usepackage{graphicx}
\usepackage{sidecap}
\usepackage{multirow}
\usepackage{booktabs}
\usepackage{bm}
\usepackage{fullpage}

\usepackage{tabularx}
\usepackage[font=small]{caption}
\usepackage{subfig}
\usepackage[linesnumbered,ruled,vlined]{algorithm2e}
\SetKwInput{KwInput}{Input}
\SetKwInput{KwOutput}{Output}

\usepackage{longtable}
\usepackage{rotating}

\usepackage{lineno,hyperref}

\graphicspath{{figures/}}
\DeclareGraphicsExtensions{.pdf,.eps,.png,.jpg,.jpeg}

\newcommand{\bfA}{\boldsymbol{A}} 
\newcommand{\tbfA}{\hat{\bfA}} 
\newcommand{\ebfA}{\tilde{\bfA}} 
\newcommand{\bfB}{\boldsymbol{B}} 
\newcommand{\tbfB}{\hat{\bfB}} 
\newcommand{\ebfB}{\tilde{\bfB}} 
\newcommand{\bfF}{\boldsymbol{F}} 
\newcommand{\tbfF}{\hat{\bfF}} 
\newcommand{\ebfF}{\tilde{\bfF}} 
\newcommand{\bfH}{\boldsymbol{H}} 
\newcommand{\tbfH}{\hat{\bfH}} 
\newcommand{\bfL}{\boldsymbol{L}} 
\newcommand{\tbfL}{\hat{\bfL}} 
\newcommand{\bfI}{\boldsymbol{I}} 
\newcommand{\bfP}{\boldsymbol{P}} 
\newcommand{\bfQ}{\boldsymbol{Q}} 

\newcommand{\bfx}{\boldsymbol{x}} 
\newcommand{\tbfx}{\hat{\boldsymbol{x}}} 
\newcommand{\ebfx}{\tilde{\boldsymbol{x}}} 
\newcommand{\bfu}{\boldsymbol{u}} 
\newcommand{\bfX}{\boldsymbol{X}} 
\newcommand{\tbfX}{\hat{\bfX}} 
\newcommand{\bfU}{\boldsymbol{U}} 
\newcommand{\bfV}{\boldsymbol{V}} 
\newcommand{\bfv}{\boldsymbol{v}} 
\newcommand{\bfXi}{\boldsymbol{\xi}} 

\newcommand{\Mbasis}{M_b}

\newcommand{\Mtrain}{M_t}

\newcommand{\Mtest}{M_{\text{test}}}
\newcommand{\MtestIn}{M'_{\text{test}}}
\newcommand{\mutest}{\mu^{\text{test}}}
\newcommand{\bfUbasis}{\bfU^{\text{b}}}
\newcommand{\bfUtrain}{\bfU}
\newcommand{\bfXbasis}{\bfX^{\text{b}}}
\newcommand{\bfXtrain}{\bfX}
\newcommand{\bbfXtrain}{\bar{\bfX}}
\newcommand{\bbfxtrain}{\bar{\bfx}}
\newcommand{\bbfx}{\bar{\bfx}}
\newcommand{\nh}{N}
\newcommand{\nr}{n}
\newcommand{\hrho}{\hat{\rho}}
\newcommand{\bfNull}{\boldsymbol{0}}
\newcommand{\bfSigma}{\boldsymbol{\Sigma}}

\newcommand{\Dcal}{\mathcal{D}}

\newcommand{\norm}[1]{\left\lVert#1\right\rVert} 

\newenvironment{keywords}%
   {\begin{trivlist}\item[]{\bfseries\sffamily Keywords:}\ }
   {\end{trivlist}}

\ihead[]{}
\ohead[]{}
\ifoot[]{}
\ofoot[]{}

\theoremstyle{definition}

\ifpdf
\hypersetup{
  pdftitle={Physics-informed regularization and structure preservation for learning stable reduced models from data with operator inference},
  pdfauthor={Nihar Sawant, Boris Kramer, Benjamin Peherstorfer}
}
\fi

\title{Physics-informed regularization and structure preservation for learning stable reduced models from data with operator inference}

\author{Nihar Sawant\thanks{Courant Institute of Mathematical Sciences, New York University} \and Boris Kramer\thanks{Department of Mechanical and Aerospace Engineering, University of California, San Diego} \and Benjamin Peherstorfer\footnotemark[1]}

\begin{document}

\maketitle

\begin{abstract}
Operator inference learns low-dimensional dynamical-system models with polynomial nonlinear terms from trajectories of high-dimensional physical systems (non-intrusive model reduction).
This work focuses on the large class of physical systems that can be well described by models with quadratic nonlinear terms and proposes a regularizer for operator inference that induces a stability bias onto quadratic models. The proposed regularizer is physics informed in the sense that it penalizes quadratic terms with large norms and so explicitly leverages the quadratic model form that is given by the underlying physics.
This means that the proposed approach judiciously learns from data and physical insights combined, rather than from either data or physics alone.
Additionally, a formulation of operator inference is proposed that enforces model constraints for preserving structure such as symmetry and definiteness in the linear terms. Numerical results demonstrate that models learned with operator inference and the proposed regularizer and structure preservation are accurate and stable even in cases where using no regularization or Tikhonov regularization leads to models that are unstable.
\end{abstract}


  \section{Introduction}
  With a lack of models and a deluge of data in science and engineering, methods for inferring models from data become ever more important. At the same time, it is increasingly recognized that relying on data alone is insufficient to learn accurate, interpretable, and robust models of science and engineering systems. Instead, a combination of data and physical insights is necessary for learning predictive models \cite{doi:10.1098/rsta.2016.0153,Willcox2021}, which has led to a surge of interest in physics-informed machine learning and scientific machine learning; see, e.g., \cite{Brunton3932,SWISCHUK2019704, raissi2019physics,doi:10.1146/annurev-fluid-010518-040547,QIAN2020132401}. In this spirit, we propose a learning method that infers low-dimensional dynamical-system models from data and induces a stability bias via a regularizer that explicitly exploits the quadratic model form given by the underlying physics.

  There is a large body of literature on learning dynamical-system models from data. We only review those that are closest to our work. First, there is system identification that originated in the systems and control community \cite{sysId}. Antoulas and collaborators introduced the Loewner approach \cite{antoulas1986scalar, mayo2007framework,interpbook}, which has been extended from linear time-invariant systems to parameterized \cite{ionita2014data}, switched \cite{doi:10.1137/17M1120233}, structured \cite{SCHULZE2018250}, delayed \cite{SCHULZE2016125}, bilinear \cite{antoulas2016model}, quadratic bilinear \cite{gosea2018data}, and polynomial \cite{doi:10.1137/19M1259171} systems as well as to learning from time-domain data \cite{peherstorfer2017data,karachalios2020bilinear}. There is also dynamic mode decomposition (DMD) \cite{schmid2010dynamic, rowley2009spectral,tu2013dynamic,kutz2016dynamic} that best-fits linear operators to state trajectories in $L_2$ norm. Methods based on Koopman operatores have been developed to extend DMD to nonlinear systems \cite{mezic2005spectral,williams2015data,brunton2016koopman}. Finally, there are sparse identification methods such as SINDy \cite{Brunton3932} and the works \cite{Schaeffer6634,doi:10.1098/rspa.2016.0446,doi:10.1137/16M1086637}. The authors of \cite{kaptanoglu2021promoting} develop a stability regularizer for SINDy that focuses on quadratic models and is motivated by Lyapunov theory. Similarly, the work \cite{erichson2019physics} adds a loss term to encourage stability of an equilibrium and so learns deep-network models that show stable behavior. Closure modeling is another research direction that recently has seen a surge of interest in data-driven methods \cite{doi:10.1146/annurev-fluid-010518-040547,doi:10.1137/16M1106419,doi:10.2514/6.2015-1287,maulik_san_2017,maulik_san_rasheed_vedula_2019} and where stabilization plays an important role \cite{BBSK17stabilizationROMextremumSeeking,doi:10.1137/17M1145136,doi:10.1137/18M1177263}.

  Our goal is to learn low-dimensional quadratic dynamical-system models and to penalize unstable models as well as preserve structure and invariances of the dynamical systems from which data are sampled. We build on operator inference \cite{peherstorfer2016data} that infers reduced models with polynomial nonlinear terms from snapshots data. Operator inference comes with recovery guarantees under certain assumptions \cite{peherstorfer2019sampling,ErrorOpInf} and it is a building block of more general learning methods that go far beyond polynomial nonlinear terms and exploit additional physical insights \cite{QKMW_2019_transform_and_learn,QIAN2020132401,BENNER2020113433}. In \cite{doi:10.2514/1.J058943}, operator inference is used together with a physics-informed lifting approach to learn a model of a large-scale combustion system, where it has been shown that regularization is important for obtaining stable models. A Tikhonov regularizer is proposed in \cite{doi:10.2514/1.J058943}, which has been further investigated in, e.g., \cite{mcquarrie2020data,ElizabethsThesis}. In contrast, we propose a regularizer that goes beyond Tikohnov regularization and that is explicitly motivated by the nature of the quadratic model form, which in turn is given by the underlying physics. Building on the insights of \cite{411100,CHESI2007326,K2020_stability_domains_QBROMs}, we penalize quadratic terms with large norms, which critically influences the stability radius of the learned models and which is also in agreement with the findings in, e.g., \cite{benner2015two,genesio1989stability}. We present numerical results that demonstrate improved stability of models learned with the proposed regularization compared to no regularization and Tikhonov regularization.

  Section~\ref{sec:Prelim} briefly describes learning low-dimensional models with operator inference and motivates this work with a synthetic example. Section~\ref{sec:PIROpInf} proposes the physics-informed regularizer and structure preservation for operator inference. The computational procedure is discussed in Section~\ref{sec:CompProc}. Numerical results in Section~\ref{sec:NumExp} demonstrate that operator inference with the proposed regularizer learns stable models even when Tikhonov regularization and models learned without regularization are unstable. Concluding remarks are in Section~\ref{sec:Conc}.

  \section{Non-intrusive model reduction with operator inference}
  \label{sec:Prelim}
  Section~\ref{section:fom} introduces the dynamical systems of interest and Section~\ref{section:fomData} discusses sampling high-dimensional state trajectories. Classical, intrusive model reduction \cite{antoulas2005approximation,rozza2007reduced,antoulas2010interpolatory,benner2015survey} requires the availability of a model of the high-dimensional dynamical system to construct a reduced model and is recapitulated in Section~\ref{section:irom}. In Section~\ref{section:nrom}, non-intrusive model reduction with operator inference is summarized, which learns reduced models from state trajectories. The problem formulation and a motivating example are given in Section~\ref{section:toy}.

  \subsection{Dynamical system with high-dimensional states}
  \label{section:fom}
  Consider a parametrized dynamical system with quadratic nonlinear terms
  \begin{equation} \label{eq:fom}
      \frac{\mathrm d}{\mathrm dt}\bfx(t; \mu) = \bfA (\mu) \bfx(t; \mu) + \bfB(\mu) \bfu(t; \mu) + \bfF (\mu) \bfx(t; \mu)^2\,,
  \end{equation}
  where $\bfA(\mu) \in \mathbb{R}^{N \times N}$ and $\bfF(\mu) \in \mathbb{R}^{N \times {N(N+1)/2}}$ are the linear and nonlinear operators, respectively. There are $p \in \mathbb{N}$ inputs that enter linearly via the input matrix $\bfB(\mu) \in \mathbb{R}^{N \times p}$. The system operators depend on a parameter $\mu \in \Dcal$ that is independent of time. The state dimension is $N \in \mathbb{N}$ and the state at time $t \in [0, T]$ is $\bfx(t; \mu) \in \mathbb{R}^N$. The $p$-dimensional input at time $t$ is $\bfu(t; \mu) \in \mathbb{R}^p$. The initial condition is denoted as $\bfx_0(\mu) \in \mathbb{R}^N$. To each state $\bfx(t; \mu) = [x_1(t; \mu),\;\dots,\;x_N(t; \mu)]^T$, there corresponds a vector $\bfx(t; \mu)^2 \in \mathbb{R}^{N(N+1)/2}$ defined as \begin{equation}
      \bfx(t; \mu)^2 = [
         \bfx^{(1)}(t; \mu)^T, \dots,
         \bfx^{(N)}(t; \mu)^T]^T
         \label{eq:Prelim:xSquare}
         \end{equation}
  where $\bfx^{(i)}(t; \mu) = x_i(t; \mu)[
         x_1(t; \mu), \dots,
         x_i(t; \mu)]^T$ for $i = 1, \dots, N$.
  The vector $\bfx(t; \mu)^2$ contains all pairwise products of components of the state vector $\bfx(t; \mu)$ up to duplicates; see, e.g., \cite{peherstorfer2016data}.

  \subsection{Collecting data}
  \label{section:fomData}
  Discretize the time domain into $0 = t_0 < t_1 < \cdots < t_K = T$ and consider the state trajectory
  \[
  \bfX(\mu) = [\bfx_0(\mu), \cdots, \bfx_K(\mu)] \in \mathbb{R}^{N \times K+1}
  \]
  for a given parameter $\mu \in \Dcal$, initial condition $\bfx_0(\mu)$, and for a given input trajectory
  \[
  \bfU(\mu) = [\bfu(t_1; \mu), \cdots, \bfu(t_K; \mu)] \in \mathbb{R}^{p \times K}\,.
  \]
  For example, a trajectory $\bfX(\mu)$ can be obtained by numerically integrating the model \eqref{eq:fom} of a dynamical system in time.
  The trajectory
  \[
  \bfX^2(\mu) = [\bfx^2_0(\mu), \cdots, \bfx^2_K(\mu)] \in \mathbb{R}^{{N(N+1)/2} \times K + 1}
  \]
  can be generated from the state trajectory $\bfX(\mu)$ following the definition \eqref{eq:Prelim:xSquare}.

  \subsection{Classical, intrusive model reduction}
  \label{section:irom}
  For each parameter $\mu \in \{ \mu_1, \dots, \mu_M \}$ in a set of $M \in \mathbb{N}$ parameters, consider $\Mbasis \in \mathbb{N}$ input trajectories $\bfUbasis_1(\mu), \dots, \bfUbasis_{\Mbasis}(\mu)$, initial conditions $\bfx_{1,0}^{\text{b}}(\mu), \dots, \bfx_{\Mbasis,0}^{\text{b}}(\mu)$, and the corresponding state trajectories $\bfXbasis_1(\mu), \dots, \bfXbasis_{\Mbasis}(\mu)$. The state trajectories for all parameters are concatenated into the snapshot matrix
  \begin{equation}
      \bfXbasis = [\bfXbasis_1(\mu_1), \dots, \bfXbasis_{\Mbasis}(\mu_1),  \dots,  \bfXbasis_1(\mu_M), \dots, \bfXbasis_{\Mbasis}(\mu_M)] \in \mathbb{R}^{\nh \times KM\Mbasis}\,.
      \label{eq:Prelim:SnapshotMatrix}
  \end{equation}
  A proper orthogonal decomposition (POD) basis of dimension $\nr \ll \nh$ is constructed from the snapshot matrix $\bfXbasis$. The basis vectors are the columns of the matrix
  \begin{equation}
  \bfV = [\bfv_1, \dots, \bfv_\nr] \in \mathbb{R}^{\nh \times \nr}\,.
  \label{eq:Prelim:PODBasis}
  \end{equation}

  To construct a projection-based reduced model via Galerkin projection, the reduced operators $\ebfA(\mu) \in \mathbb{R}^{\nr \times \nr}$ and $\ebfB(\mu) \in \mathbb{R}^{\nr \times p}$ are obtained via
  \begin{equation}
  \ebfA(\mu) = \bfV^T \bfA(\mu) \bfV\,,\qquad \ebfB = \bfV^T \bfB(\mu)\,.
  \label{eq:Prelim:ReducedOps}
  \end{equation}
  The reduced quadratic operator $\ebfF(\mu) \in \mathbb{R}^{\nr \times \nr(\nr+1)/2}$ is constructed in a similar fashion via projection as described in, e.g., \cite{peherstorfer2019sampling}. Thus, for each parameter $\mu \in \{\mu_1,\;,\dots,\;,\mu_M \}$, one obtains a reduced model
  \begin{equation}
      \frac{\mathrm{d}}{\mathrm dt}\ebfx(t, \mu) = \ebfA(\mu) {\ebfx}(t; \mu) + \ebfB (\mu) \bfu(t; \mu) + \ebfF(\mu) {\ebfx}^2(t; \mu)\,,
      \label{eq:Prelim:ROM}
  \end{equation}
  with reduced state $\ebfx(t; \mu) \in \mathbb{R}^{\nr}$. We refer to \cite{antoulas2005approximation,rozza2007reduced,benner2015survey} for more details about classical, intrusive model reduction.

  \subsection{Learning low-dimensional models from data with operator inference}
  \label{section:nrom}
  Constructing the reduced operators $\ebfA(\mu), \ebfB(\mu), \ebfF(\mu)$ via \eqref{eq:Prelim:ReducedOps} is an intrusive process because it requires access to high-dimensional operators $\bfA(\mu), \bfB(\mu)$, $\bfF(\mu)$ in either implicit or explicit form. In contrast, operator inference \cite{peherstorfer2016data} aims to learn reduced operators from trajectories of the dynamical system \eqref{eq:fom}, without requiring access to the high-dimensional operators.

  Let $\bfV$ be a basis matrix. Notice that such a basis matrix $\bfV$ can be constructed purely from the snapshot matrix \eqref{eq:Prelim:SnapshotMatrix} in many situations, without having available the high-dimensional operators. Then, the intrusive projection step \eqref{eq:Prelim:ReducedOps} is replaced with a non-intrusive least-squares regression problem. First, for each training parameter $\mu \in \{\mu_1, \dots, \mu_M\}$, consider $\Mtrain$ training input trajectories $\bfUtrain_1(\mu), \dots, \bfUtrain_{\Mtrain}(\mu)$ with initial conditions $\bfx_{1,0}(\mu), \dots, \bfx_{\Mtrain,0}(\mu)$ and the corresponding training state trajectories $\bfXtrain_1(\mu), \dots, \bfXtrain_{\Mtrain}(\mu)$. Second, the training trajectories are projected onto the reduced space via
  \begin{equation*}
      \bbfXtrain_i(\mu) = \bfV^T \bfXtrain_i(\mu)
  \end{equation*}
  to obtain the projected training trajectories $\bbfXtrain_i(\mu) = [\bbfxtrain_{i,1}(\mu), \dots, \bbfxtrain_{i,K}(\mu)] \in \mathbb{R}^{\nr \times K}$ for $i = 1, \dots, \Mtrain$. Third, the operators $\tbfA(\mu), \tbfB(\mu)$ and $\tbfF(\mu)$ are fitted via least-squares regression to the projected training trajectories
  \begin{equation*}
      \min_{\tbfA(\mu), \tbfB(\mu), \tbfF(\mu)} J(\tbfA(\mu), \tbfB(\mu), \tbfF(\mu))
      \end{equation*}
      with objective function
  \begin{equation}
      J(\tbfA(\mu), \tbfB(\mu), \tbfF(\mu)) = \sum_{i = 1}^{\Mtrain}\sum_{k=1}^K \norm{\bbfx_{i,k}^{\prime}(\mu) - \tbfA(\mu) \bbfx_{i,k}(\mu) - \tbfB(\mu) \bfu_{i,k}(\mu) - \tbfF(\mu) \bbfx_{i,k}^2(\mu)}_2^2\,,
      \label{eq:Prelim:OpInfJ}
  \end{equation}
  where $\bfu_{i,k}(\mu)$ is the input at time step $k$ of the $i$th training trajectory $\bfU_i = [\bfu_{i,1}(\mu), \dots, \bfu_{i,K}]$ for $i = 1, \dots, \Mtrain$. The quantity $\bbfx_{i,k}^{\prime} \in \mathbb{R}^{\nr}$ denotes a numerical approximation of the time derivative of the projected state at time $k$ of the $i$th trajectory, such as a first-order finite difference approximation
  \begin{equation}
  \bbfx_{i,k}^{\prime} = \frac{\bbfx_{i,k} - \bbfx_{i,k-1}}{\delta t}\,,
  \label{eq:Prelim:OpInfDiffApprox}
  \end{equation}
  with time-step size $\delta t > 0$. The inferred operators $\tbfA(\mu), \tbfB(\mu), \tbfF(\mu)$ are then used to assemble a low-dimensional model
  \begin{equation}
      \frac{\mathrm d}{\mathrm dt} \tbfx(t; \mu) =  \tbfA(\mu) \tbfx_k(t; \mu) + \tbfB(\mu) \bfu_k(\mu) + \tbfF(\mu) \tbfx_k^2(t; \mu),
      \label{eq:Prelim:OpInfROM}
  \end{equation}
  where $\tbfx_k(t; \mu) \in \mathbb{R}^{\nr}$ is the state at time $t$.

  The operator inference process is repeated for each parameter in the training set $\{\mu_1, \dots, \mu_M\}$ to compute the corresponding inferred operators. For a new parameter $\mu \in \Dcal \setminus \{\mu_1, \dots, \mu_M\}$, the operators $\tbfA(\mu), \tbfB(\mu), \tbfF(\mu)$ are obtained via interpolation. We refer to \cite{peherstorfer2016data} for details and to \cite{lin2019riemannian} for interpolating between reduced operators in model reduction in general.

  \subsection{Stability of inferred models}
  \label{section:toy}
  We demonstrate operator inference on a toy example. Consider the dynamical system
  \begin{equation}
  \frac{\mathrm d}{\mathrm dt} \bfx(t) = \bfA(\mu)\bfx(t) + \bfB u(t) + \bfF \bfx^2(t)\,,
  \label{eq:Prelim:ProbForm:System}
  \end{equation}
  where $\bfB \in \mathbb{R}^{\nh \times 1}$ and $\bfF \in \mathbb{R}^{\nh \times \nh(\nh + 1)/2}$ have entries that are realizations of the uniform distribution in $[0, 1]$. The dimension is $\nh = 128$. The linear operator in \eqref{eq:Prelim:ProbForm:System} is $\bfA(\mu) = -\mu (\bfA_s + \bfA_s^T + 2\nh \bfI)$, where $\bfI$ is the identity matrix and $\bfA_s \in \mathbb{R}^{\nh \times \nh}$ is a matrix that has as entries realizations of the uniform distribution in $[0, 1]$. The matrix $\bfA(\mu)$ is symmetric negative definite with probability 1.
  The parameter domain is $\Dcal = [0.1, 1]$ and
  end time is $T = 1$. We discretize \eqref{eq:Prelim:ProbForm:System} with time-step size $\delta t = 10^{-3}$ and explicit Euler.
  For each training parameter $\mu \in \{0.1, \; \dots, \; 1.0\}$, we generate a single ($\Mbasis = 1$) input trajectory $\bfUbasis_1(\mu)$, whose entries are random with a uniform distribution in $[0, 2]$, and an initial condition $\bfx^{\text{b}}_{1,0}$, whose entries follow a uniform distribution in $[0, 1]$. The corresponding state trajectories are $\bfXbasis_1(\mu_1), \dots, \bfXbasis_{1}(\mu_M)$.
  A basis matrix $\bfV \in \mathbb{R}^{\nh \times \nr}$ is then constructed from the corresponding snapshots as described in Section~\ref{section:irom}. The reduced basis is generated for dimension $\nr = 2, \; \dots, \; 10$.
  For parameter $\mu = 0.7$, we then construct $\Mtrain = 3$ training inputs $\bfUtrain_1, \dots, \bfUtrain_{\Mtrain}$ with training initial conditions $\bfx_{1,0}, \dots, \bfx_{\Mtrain,0}$, which are sampled from the same distributions as the inputs and initial conditions for the basis construction. The corresponding training state trajectories are $\bfXtrain_1, \dots, \bfXtrain_{\Mtrain}$, to which we apply operator inference as described in Section~\ref{section:nrom}. We use a first-order forward difference scheme to approximate the time derivative as in \eqref{eq:Prelim:OpInfDiffApprox}. Additionally, for benchmarking purposes, we also construct a reduced model \eqref{eq:Prelim:ROM} via the intrusive process described in Section~\ref{section:irom}.

\begin{figure}
    \centering
    \subfloat[training error \eqref{eq:Prelim:SynTrainError}]{{\LARGE\resizebox{0.48\columnwidth}{!}{\input{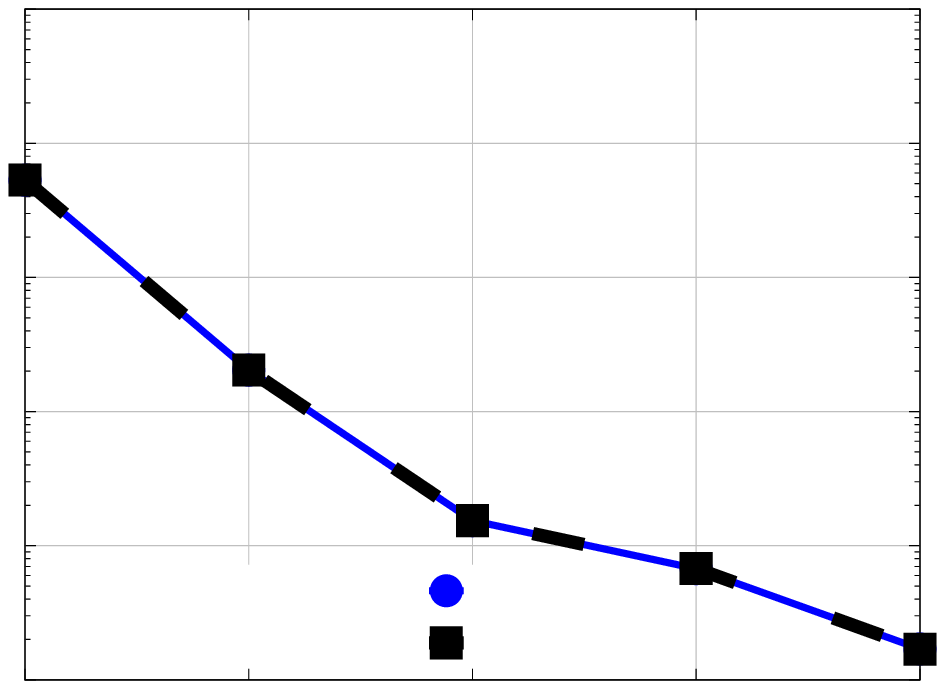}}}
    \label{fig:sample-train-noreg}}
    \hfill
    \subfloat[test error \eqref{eq:Prelim:SynTestError}]{{\LARGE\resizebox{0.48\columnwidth}{!}{\input{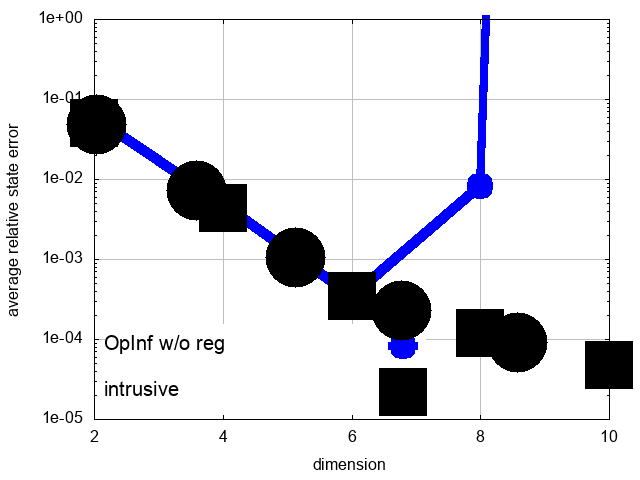}}}
    \label{fig:sample-test-noreg}}
    \caption{Synthetic example: (a) The training error for the model learned via operator inference (OpInf) matches the error of the model obtained with intrusive model reduction. (b) When tested for a different input and initial condition at the same parameter, the model learned via operator inference (without regularization) is inaccurate and unstable in this example.}
  \end{figure}

  Figure~\ref{fig:sample-train-noreg} shows the training error of the operator-inference model
  \begin{equation}
  e_{\text{train}} = \sum_{i=1}^{\Mtrain} \frac{\| \bfV \tbfX_i - \bfX_i \|_F}{\| \bfX_i \|_F}\,,
  \label{eq:Prelim:SynTrainError}
  \end{equation}
  which indicates that the operator-inference model achieves a comparable error decay as the reduced model obtained from intrusive model reduction. However, if we simulate the operator-inference model at a test input trajectory $\bfU^{\text{test}}$, whose entries are sample uniformly in $[0, 10]$, and test initial condition $\bfx_0^{\text{test}}$, with entries sampled uniformly in $[0, 1]$, and plot the error
  \begin{equation}
  e_{\text{test}} = \frac{\| \bfV \tbfX^{\text{test}} - \bfX^{\text{test}} \|_F}{\| \bfX^{\text{test}} \|_F}\,,
  \label{eq:Prelim:SynTestError}
  \end{equation}
  in Figure~\ref{fig:sample-test-noreg}, then an instability can be observed, compared to the reduced model from intrusive model reduction. The results indicate that operator inference is prone to overfitting, which can result in unstable behavior at test inputs as in this motivating example.

  \section{A physics-informed regularizer for operator inference}
  \label{sec:PIROpInf}
  We propose a physics-informed regularizer for operator inference that penalizes unstable dynamical-system models. In Section~\ref{sec:PIOpInf:Stability}, we recapitulate the definition of the stability radius of quadratic dynamical-system models with Lyapunov stability criteria. We then propose a regularizer that penalizes models with small stability radii in Section~\ref{sec:PIOpInf:Reg}. Additionally, in Section~\ref{sec:PIOpInf:Structure}, we propose to combine the phyics-informed regularizer with structure preservation. In the following, for each parameter $\mu \in \{\mu_1, \dots, \mu_M\}$ in the training set, an operator-inference model is learned separately as discussed in Section~\ref{section:nrom}. Thus, in this section, the parameter dependence of quantities is not explicitly denoted.

  \subsection{Stability radius of quadratic models of dynamical systems}
  \label{sec:PIOpInf:Stability}
  We closely follow \cite{K2020_stability_domains_QBROMs} to define the concept of stability radius for quadratic models of dynamical systems. In fact, the work \cite{K2020_stability_domains_QBROMs} is the motivation for the proposed physics-informed regularizer.

  \subsubsection{Stability domain of quadratic models}
  Consider the autonomous quadratic reduced model
  \begin{equation}
  \frac{\mathrm d}{\mathrm dt} \tbfx(t) = \tbfA \tbfx(t) + \tbfH (\tbfx(t) \otimes \tbfx(t))\,
  \label{eq:PIOpInf:AutoQSystem}
  \end{equation}
  where $\otimes$ denotes the Kronecker product. Notice that \eqref{eq:PIOpInf:AutoQSystem} is a different representation of an autonomous version of
  the quadratic reduced model \eqref{eq:Prelim:ROM} and the quadratic operator-inference model \eqref{eq:Prelim:OpInfROM} if there is no input.
  Similarly, we can represent the quadratic full model defined in \eqref{eq:fom} in the form \eqref{eq:PIOpInf:AutoQSystem} if there is no input.

  Let now without loss of generality $\tbfx_e = \bfNull$ be an equilibrium point of \eqref{eq:PIOpInf:AutoQSystem}, i.e., $\tbfA \tbfx_e + \tbfH( \tbfx_e \otimes \tbfx_e) = \bfNull$. The domain of attraction $\mathcal{A}(\tbfx_e)$ of the equilibrium $\tbfx_e$ is then defined as the set of initial conditions that lead to the equilibrium point $\tbfx_e$ as a steady state, i.e.,
  \begin{equation*}
      \mathcal{A}(\tbfx_e) = \{\tbfx_0 : \lim_{t \to \infty} \tbfx(t) = \tbfx_e\}\,,
  \end{equation*}
  where $\tbfx(t)$ is the state at time $t$ of \eqref{eq:PIOpInf:AutoQSystem} with initial condition $\tbfx_0$. Directly working with the stability domain is challenging from an analytic and computational point of view and thus one typically resorts to deriving subsets $D \subseteq \mathcal{A}(\tbfx_e)$. To measure a subset $D$, we build on the Lyapunov theory to derive a stability radius.

  If there exists a Lyapunov function $\nu: \mathbb{R}^{\nh} \to \mathbb{R}^+$ that is continuously differential and that satisfies
  \begin{equation*}
      \nu(\tbfx) > 0, \; \dot{\nu}(\tbfx) < 0, \quad \forall \tbfx \in \mathcal{A}(\tbfx_e)\,,
  \end{equation*}
  then model \eqref{eq:PIOpInf:AutoQSystem} is locally asymptotically stable about $\tbfx_e$.
  Here, $\dot{\nu}(\tbfx)$ means $\dot{\nu}(\tbfx) = \frac{\mathrm d\nu}{\mathrm d\tbfx} \boldsymbol{\hat{f}}(\tbfx)$, where $\boldsymbol{\hat{f}}(\tbfx) = \tbfA\tbfx + \tbfH(\tbfx \otimes \tbfx)$ is the right-hand side function of the corresponding dynamical system.
  As shown in \cite{411100,CHESI2007326,K2020_stability_domains_QBROMs}, given a Lyapunov function $\nu$, an estimate $D(\rho) \subseteq \mathcal{A}(\tbfx_e)$ of the domain of attraction $\mathcal{A}(\tbfx_e)$ is given by
      \begin{equation*}
          D(\rho) = \{ \tbfx: \nu(\tbfx) \le \rho^2, \; \dot{\nu}(\tbfx) < 0 \}\,,
      \end{equation*}
  where we refer to $\rho$ as the stability radius.

  \subsubsection{Estimating stability radius}
  Consider an autonomous quadratic model \eqref{eq:PIOpInf:AutoQSystem} with Lyapunov function $\nu(\tbfx) = \tbfx^T \bfP \tbfx$, where $\bfP \in \mathbb{R}^{n \times n}$ is a symmetric positive definite matrix that satisfies
  \begin{equation}
  \bfL\bfL^T = -\tbfA^T \bfP - \bfP \tbfA\,,
  \label{eq:PIOpInf:LyapunovEq}
  \end{equation}
  for an arbitrary matrix $\bfL \in \mathbb{R}^{n \times n}$.
  The derivative of the Lyapunov function along a trajectory is
  \[
  \dot{\nu}(\tbfx) = \dot{\tbfx}^T \bfP \tbfx + \tbfx^T \bfP \dot{\tbfx}\,.
  \]
  Building on \cite[Proposition~3.1]{K2020_stability_domains_QBROMs}, we obtain the radius
  \begin{equation}
  \hrho = \frac{\sigma_{\text{min}}(\bfL)}{2\sqrt{\|\bfP\|_F}\| \tbfH\|_F}
  \label{eq:PIOpInf:Radius}
  \end{equation}
  and that $D(\hrho) \subseteq \mathcal{A}(\tbfx_e)$ is a subset of $\mathcal{A}(\tbfx_e)$, if $\tbfA$ is Hurwitz, i.e., the real parts of all eigenvalues of $\tbfA$ are negative.
  Notice that in contrast to the 2-norm $\|\cdot\|_2$ used in \cite[Proposition~3.1]{K2020_stability_domains_QBROMs}, we state the radius \eqref{eq:PIOpInf:Radius} with respect to the Frobenius norm $\|\cdot\|_F$, which leads to an operator inference problem that can be solved more efficiently than when working with the $\|\cdot\|_2$ norm.

  \subsection{Operator inference with physics-informed regularizer}
  \label{sec:PIOpInf:Reg}
  The stability radius $\hrho$, which is derived in \eqref{eq:PIOpInf:Radius}, grows inversely proportional to the norm $\|\tbfH\|_F$ of the quadratic term $\tbfH \in \mathbb{R}^{\nr \times \nr^2}$ in \eqref{eq:PIOpInf:AutoQSystem}. Notice that other results on stability analysis for quadratic systems, e.g.,~\cite{benner2015two,genesio1989stability}, also show that a small norm of the quadratic term can increase the stability radius.

  The models that we infer with operator inference have the form \eqref{eq:Prelim:OpInfROM} and thus the quadratic term $\tbfF$ is of dimension $\nr \times \nr(\nr + 1)/2$. However, models of the form \eqref{eq:PIOpInf:AutoQSystem} with $\tbfF$ can be transformed into models with quadratic terms $\tbfH (\tbfx(t) \otimes \tbfx(t))$ such that $\tbfF \tbfx^2 = \tbfH (\tbfx \otimes \tbfx)$ holds for any $\tbfx \in \mathbb{R}^{\nr}$ and $\|\tbfH\|_F \leq \|\tbfF\|_F$ holds as well. We can construct $\tbfH$, such that $\|\tbfH\|_F = \|\tbfF\|_F$, by filling the additional columns of $\tbfH$ with zeros.
  Thus, we obtain that if we regularize the norm $\|\tbfF\|_F$, we also regularize the norm $\|\tbfH\|_F$ of a corresponding $\tbfH$, which in turn means that the denominator of the radius $\hrho$ is regularized. This leads to the optimization problem for inferring model \eqref{eq:Prelim:OpInfROM} with operator inference and the proposed physics-informed regularizer (PIR-OpInf)
  \begin{equation}
      \min_{\hat{\bfA}, \hat{\bfB}, \hat{\bfF}} J(\hat{\bfA},\hat{\bfB},\hat{\bfF}, \lambda) + \lambda \|\hat{\bfF}\|_F^2\,,
      \label{eq:PIOpInf:ObjPIR}
  \end{equation}
  with $J$ defined in \eqref{eq:Prelim:OpInfJ} and $\lambda > 0 $ being a regularization parameter. Notice that increasing $\lambda$ means more severely penalizing the norm $\|\tbfF\|_F$, which in turn leads to a potential increase of the radius $\hrho$ and thus a more stable inferred model in the sense of Lyapunov.

  The PIR-OpInf problem \eqref{eq:Prelim:OpInfJ} imposes no constraints on the linear operator $\hat{\bfA}$. In particular, there is no guarantee that the inferred $\hat{\bfA}$ is Hurwitz and thus there can exist eigenvalues with non-negative real parts. To ensure a linear operator that is Hurwitz, we apply an eigenvalue reflection as a post-processing step. Let $\hat{\bfA} = \bfQ_A\bfSigma_A\bfQ_A^{-1}$ be the eigendecomposition of $\hat{\bfA}$. If $\hat{\bfA}$ is not diagonalizable, we reduce the dimension $\nr$ until a matrix $\hat{\bfA}$ is inferred that is diagonalizable. Notice that this process stops in a finite number of steps because $\hat{\bfA}$ is diagonalizable if $\nr = 1$ and $\hat{\bfA}$ is non-zero. Without loss of generality, let $\sigma_1, \dots, \sigma_r$ be eigenvalues with non-negative real parts and let $\sigma_{r+1}, \dots, \sigma_n$ be all other eigenvalues. Denote with $\mathfrak{R}(\sigma)$ and $\mathfrak{I}(\sigma)$ the real and imaginary part, respectively, for a complex number $\sigma \in \mathbb{C}$. Then, we replace $\hat{\bfA}$ with the matrix
  \[
  \bfQ_A \operatorname{diag}(-\epsilon + \mathfrak{I}(\sigma_1), \dots, -\epsilon + \mathfrak{I}(\sigma_r), \sigma_{r + 1}, \dots, \sigma_{\nr})\bfQ_A^{-1}\,,
  \]
  which replaces the positive real parts of the eigenvalues with a negative real number given by the small positive threshold $\epsilon > 0$. Notice that other post-processing strategies can be applied to obtain a Hurwitz linear operator; we refer to \cite{HIGHAM1988103,KALASHNIKOVA2014251,7496283}. Note further that the post-processing also needs to be applied after interpolating at a new parameter $\mu \in \Dcal$ outside of the training set; cf.~Section~\ref{section:nrom}.

  \subsection{Operator inference with structure preservation}
  \label{sec:PIOpInf:Structure}
  Structure can be imposed on the linear operator by adding hard constraints to the operator inference problem. We focus on problems that lead to symmetric negative definite linear operators and thus we consider the constrained problem
  \begin{equation}
  \begin{aligned}
      \min_{\hat{\bfA}, \hat{\bfB}, \hat{\bfF}} ~~~& J(\hat{\bfA},\hat{\bfB},\hat{\bfF}) + \lambda \|\hat{\bfF}\|_F^2\,,\\
      \text{such that } ~~~& \hat{\bfA} - \epsilon\bfI \preceq 0\,,
      \end{aligned}
      \label{eq:PIOpInf:ObjSPIR}
  \end{equation}
  where $\hat{\bfA} - \epsilon \bfI \preceq 0$ means that $\hat{\bfA} - \epsilon \bfI$ is symmetric negative semi-definite. The matrix $\bfI$ is the identity and $\epsilon > 0$ is a margin that guarantees that $\hat{\bfA}$ is definite, rather than semi-definite. We refer to \eqref{eq:PIOpInf:ObjSPIR} as the SPIR-OpInf problem, where the S stands for ``structure.'' Problem \eqref{eq:PIOpInf:ObjSPIR} is a semi-definite program, for which efficient algorithms exist \cite{Boyd}.

  Instead of imposing a margin $\epsilon$ to guarantee definiteness of $\hat{\bfA}$ in \eqref{eq:PIOpInf:ObjSPIR}, one can solve problem \eqref{eq:PIOpInf:ObjSPIR} with the constraint $\hat{\bfA} \preceq 0$ and subsequently apply an analogous post-processing step as in Section~\ref{sec:PIOpInf:Reg}. Because symmetry is enforced by $\hat{\bfA} \preceq 0$, it is guaranteed that $\hat{\bfA}$ is diagonalizable. The post-processing described in Section~\ref{sec:PIOpInf:Reg} preserves symmetry and thus the result is a symmetric negative definite matrix after the post-processing; see also \cite{HIGHAM1988103}. However, in the following, we will impose a margin $\epsilon$ and therefore do not need a post-processing step.

  Other structures in the linear operator can be preserved in an analogous way. For example, another common structure is skew-symmetry of $\hat{\bfA}$ which can be formulated as a linear constraint; we leave such other constraints to future work.
  Recall that it is required to interpolate between inferred operators when a model at a parameter $\mu$ outside of the training set is required; cf.~Section~\ref{section:nrom}. In case of structure-preserving operator inference, the corresponding operator interpolation schemes have to preserve the operator structure. We discuss such an interpolation scheme for symmetric negative definite matrices in Section~\ref{sec:Comp:Structure}.

  \section{Computational procedure of physics-informed operator inference}
  \label{sec:CompProc}
  In Section~\ref{sec:Comp:Structure}, we briefly recapitulate an interpolation scheme that preserve symmetric definiteness of matrices, which is critical for constructing operators at new parameters outside of the training set in SPIR-OpInf. To select a regularization parameter for PIR-OpInf \eqref{eq:PIOpInf:ObjPIR} and SPIR-OpInf \eqref{eq:PIOpInf:ObjSPIR}, we propose a parameter-selection scheme in Section~\ref{sec:Comp:CV}. Section~\ref{sec:Comp:Algo} presents Algorithm~\ref{algo:PIROpInf} that summarizes the computational procedure for operator inference with physics-informed regularization and structure preservation.

  \subsection{Interpolation of structure-preserving operator-inference models}
  \label{sec:Comp:Structure}
  In SPIR-OpInf introduced in Section~\ref{sec:PIOpInf:Structure}, the definiteness and symmetry constraints in the optimization problem \eqref{eq:PIOpInf:ObjSPIR} ensure that for each training parameter $\mu_1, \dots, \mu_M$ a model is obtained with a linear operator that is symmetric negative definite. When we interpolate the trained models at a new parameter $\mu \in \Dcal \setminus \{\mu_1, \dots, \mu_M\}$ outside of the training set, however, we have to ensure that the interpolated linear operator is symmetric negative definite as well. There are various interpolation schemes in model reduction that preserve such structure, see, e.g., \cite{degroote2010interpolation, amsallem2008interpolation}.
  We build on the Log-Cholesky averaging method presented in \cite{lin2019riemannian}.

  Given are $M$ symmetric negative definite matrices $\tbfA(\mu_1), \dots, \tbfA(\mu_M)$ at parameters $\mu_1, \dots, \mu_M$. We compute the Cholesky factors $\tbfL(\mu_i)$ such that $\tbfA(\mu_i) = -\tbfL(\mu_i) \tbfL(\mu_i)^T$ for $i = 1, \dots, M$. The Cholesky factors $\tbfL(\mu_i)$ are then split into
  \begin{equation*}
     \tbfL(\mu_i) = \lfloor \tbfL(\mu_i) \rfloor + \operatorname{diag}(\tbfL(\mu_i))\,,\qquad i = 1, \dots, M\,,
  \end{equation*}
  where $\operatorname{diag}(\tbfL(\mu_i))$ is the diagonal matrix with the same the diagonal as $\tbfL(\mu_i)$ and $\lfloor \tbfL(\mu_i) \rfloor$ is its remaining strictly lower triangular part. The interpolated matrix $\tbfA(\mu)$ at a new parameter $\mu$ is
  \begin{equation*}
      \tbfA(\mu) = -\tbfL(\mu)\tbfL(\mu)^T\,,
  \end{equation*}
  where the Cholesky factor $\tbfL(\mu)$ is obtained as
  \[
  \tbfL(\mu) = \mathcal{I}(\mu; \lfloor \tbfL(\mu_1) \rfloor, \dots, \lfloor \tbfL(\mu_M) \rfloor) + \operatorname{exp}\left(\mathcal{I}\left(\mu; \log\left(\operatorname{diag}\left(\tbfL(\mu_1)\right)\right), \dots, \log\left(\operatorname{diag}\left(\tbfL(\mu_M)\right)\right)\right)\right)\,.
  \]
  The operator $\mathcal{I}$ denotes linear interpolation of the matrix entries at $\mu$ and $\exp(\cdot)$ and $\log(\cdot)$ are the matrix exponential and logarithm, respectively.

  \subsection{A parameter-selection scheme for PIR-OpInf and SPIR-OpInf}
  \label{sec:Comp:CV}
  Let $\mu_1, \dots, \mu_M$ be the training parameters and recall that $\bfXtrain_1(\mu_i), \dots, \bfXtrain_{\Mtrain}(\mu_i)$ are the training trajectories with input trajectories $\bfUtrain_1(\mu_i), \dots, \bfUtrain_{\Mtrain}(\mu_i)$, respectively, for $i = 1, \dots, M$; cf.~Section~\ref{section:irom}.
  Define the minimum $\lambda_{\text{min}}$ and maximum $\lambda_{\text{max}}$ of the regularization parameter and discretize the interval $[\lambda_{\text{min}}, \lambda_{\text{max}}] \subset \mathbb{R}$ with $m$ points
  \begin{equation}
  \lambda_{\text{min}} = \lambda_1 < \dots < \lambda_m = \lambda_{\text{max}}.
  \label{eq:Comp:LambdaDisc}
  \end{equation}

  For each $\lambda_i$, we learn a model $\hat{\Sigma}_{ij}$ with PIR-OpInf \eqref{eq:PIOpInf:ObjPIR} for $\mu_j$, with $i = 1, \dots, m$ and $j = 1, \dots, M$. Then, for $j = 2, \dots, M-1$ and for $i = 1, \dots, m$, we derive $\hat{\Pi}_{ij}$ by interpolating between models
  \begin{equation}
  \hat{\Sigma}_{i,1}, \dots, \hat{\Sigma}_{i, j-1}, \hat{\Sigma}_{i, j+1}, \dots, \hat{\Sigma}_{i,M}
  \label{eq:Comp:Reg:AllModels}
  \end{equation}
  corresponding to parameters $\mu_1, \dots, \mu_{j - 1}, \mu_{j + 1}, \dots, \mu_M$, i.e., the parameter $\mu_j$ corresponding to model $\hat{\Sigma}_{ij}$ is left out from the interpolation process. The interpolation is structure preserving if necessary; cf.~Section~\ref{sec:Comp:Structure}. Notice that all models in \eqref{eq:Comp:Reg:AllModels} are trained with the same regularization parameter $\lambda_i$. The interpolated model $\hat{\Pi}_{ij}$ is integrated in time with the input trajectories $\bfUtrain_1(\mu_j), \dots, \bfUtrain_{\Mtrain}(\mu_j)$ corresponding to parameter $\mu_j$ to obtain the trajectories $\hat{\bfX}^{(i)}_1(\mu_j), \dots, \hat{\bfX}^{(i)}_{\Mtrain}(\mu_j)$ and the error
  \begin{equation}
  e_{ij}^{\text{val}} = \sum_{\ell = 1}^{\Mtrain} \frac{\|\bfV\hat{\bfXtrain}_{\ell}^{(i)}(\mu_j) - \bfXtrain_{\ell}(\mu_j)\|_F}{\|\bfXtrain_{\ell}(\mu_j)\|_F}
  \label{sec:Comp:ParamSelErr}
  \end{equation}
  is assigned to the pair of regularization parameter $\lambda_i$ and parameter $\mu_j$, where $\bfV$ is the basis matrix. We then pick $\lambda^*$ by solving
  \begin{equation}
  \operatorname*{arg\,min}_{i = 1, \dots, m}\,\,\, \frac{1}{M-2}\sum_{i = 2}^{M-1} e_{ij}^{\text{val}}\,.
  \label{sec:Comp:ParamSelOpti}
  \end{equation}
  The same procedure is applied in case of SPIR-OpInf \eqref{eq:PIOpInf:ObjSPIR}.

  Notice that error due to interpolating between models enters the validation error \eqref{sec:Comp:ParamSelErr} and thus the selection criterion \eqref{sec:Comp:ParamSelOpti} for $\lambda^*$. This is in contrast to other parameter-selection schemes for operator inference that are either formulated in parameter-independent settings or ignore the parameter dependency in the selection process \cite{doi:10.2514/1.J058943,mcquarrie2020data}.

  \subsection{Algorithm of operator inference with physics-informed regularizer and structure preservation}
  \label{sec:Comp:Algo}
  Algorithm~\ref{algo:PIROpInf} summarizes the computational procedure of the proposed approach. Inputs are the basis matrix $\bfV$, which is constructed from trajectories as described in Section~\ref{section:irom}, and the training trajectories $\bfXtrain_1(\mu_i), \dots, \bfXtrain_{\Mtrain}(\mu_i)$ and inputs $\bfUtrain_1(\mu_i), \dots, \bfUtrain_{\Mtrain}(\mu_i)$ for the training parameter $\mu_i$ with $i = 1, \dots, M$. In the nested for loop, models are generated with either PIR-OpInf \eqref{eq:PIOpInf:ObjPIR} or SPIR-OpInf \eqref{eq:PIOpInf:ObjSPIR} for all pairwise combinations of regularization parameters defined in \eqref{eq:Comp:LambdaDisc} and training parameters $\mu_1, \dots, \mu_M$. Then, the validation error \eqref{sec:Comp:ParamSelErr} is computed and the index $i^*$ of the regularization parameter $\lambda_{i^*}$ that minimizes the validation error is determined. The corresponding inferred models are returned.

  \begin{algorithm}[t]
  \DontPrintSemicolon
  \KwInput{basis $\bfV$, inputs $\bfUtrain_1(\mu_j), \dots, \bfUtrain_{\Mtrain}(\mu_j)$ and trajectories $\bfXtrain_1(\mu_j), \dots, \bfXtrain_{\Mtrain}(\mu_j)$ for $j = 1, \dots, M$}
  \KwOutput{inferred operators $\tbfA(\mu_j), \tbfB(\mu_j), \tbfF(\mu_j)$ for $j = 1, \dots, M$}
  \For{$i = 1, \dots, m$}{
  \For{$j = 1, \dots, M$}{
  Infer operators $\tbfA^{(i)}(\mu_j), \tbfB^{(i)}(\mu_j), \tbfF^{(i)}(\mu_j)$ with either PIR-OpInf \eqref{eq:PIOpInf:ObjPIR} or SPIR-OpInf \eqref{eq:PIOpInf:ObjSPIR} and regularization parameter $\lambda_i$ defined in \eqref{eq:Comp:LambdaDisc} and training parameter $\mu_j$\;
  }
  }
  Compute validation error \eqref{sec:Comp:ParamSelErr} for $i = 1, \dots, m$ and $j = 2, \dots, M - 1$\;
  Pick $\lambda^* = \lambda_{i^*}$ with index $i^*$ as in \eqref{sec:Comp:ParamSelOpti} that minimizes validation error \;
  Set $\tbfA(\mu_j) = \tbfA^{(i^*)}(\mu_j), \tbfB(\mu_j) = \tbfB^{(i^*)}(\mu_j), \tbfF(\mu_j) = \tbfF^{(i^*)}(\mu_j)$ for $j = 1, \dots, M$ \;
  \Return $\tbfA(\mu_1), \dots, \tbfA(\mu_M), \tbfB(\mu_1), \dots, \tbfB(\mu_M), \tbfF(\mu_1), \dots, \tbfF(\mu_M)$
  \caption{Operator inference with physics-informed regularizer and structure preservation}
  \label{algo:PIROpInf}
  \end{algorithm}

  \section{Numerical experiments}
  \label{sec:NumExp}
  In this section, we compare operator inference with the proposed physics-informed regularizer (PIR-OpInf) and structure preservation (SPIR-OpInf) to Tikhonov regularization and operator inference without regularization. Section~\ref{sec:NumExp:Syn} revisits the synthetic example from Section~\ref{section:toy}. Section~\ref{sec:NumExp:Burgers} and Section~\ref{sec:NumExp:React} show experiments with the Burgers' equation and a reaction-diffusion problem in a pipe. The proposed approach depends on a small, positive threshold $\epsilon > 0$, e.g., for the post-processing in PIR-OpInf (cf.~Section~\ref{sec:PIOpInf:Reg}) and for the margin in SPIR-OpInf \eqref{eq:PIOpInf:ObjSPIR}, which we set to $\epsilon = 10^{-10}$ in all of the following experiments.

  \subsection{Synthetic example}
  \label{sec:NumExp:Syn}
  Consider again the synthetic example introduced in Section~\ref{section:nrom}. We now apply PIR-OpInf with the parameter-selection procedure discussed in Section~\ref{sec:Comp:CV}. For each dimension $\nr \in \{2, 4, 6, 8, 10\}$, we sweep over $m = 51$ regularization parameters that are log-uniformly distributed in the interval $[10^{-15}, 10^5]$. The selected regularization parameters are $\lambda^* = 10^{-10}, 1.58 \times 10^{-7}, 10^{-8}, 3.98 \times 10^{-9}, 1.58 \times 10^{-9}$ for dimensions $\nr = 2, 4, 6, 8, 10$, respectively.
  We choose the test parameter set $\{\mutest_1, \dots, \mutest_{\Mtest}\}$ of $\Mtest = 7$ test parameters that are equidistantly chosen in $\Dcal$, where for each test parameter a test input trajectory is constructed with entries sampled uniformly in $[0, 10]$ and a test initial condition with entries sampled uniformly in $[0, 1]$, cf.~Section~\ref{section:toy}.

\begin{figure}
    \centering
    \subfloat[test error \eqref{eq:NumExp:Syn:TestError}]{{\Large\resizebox{0.48\columnwidth}{!}{\input{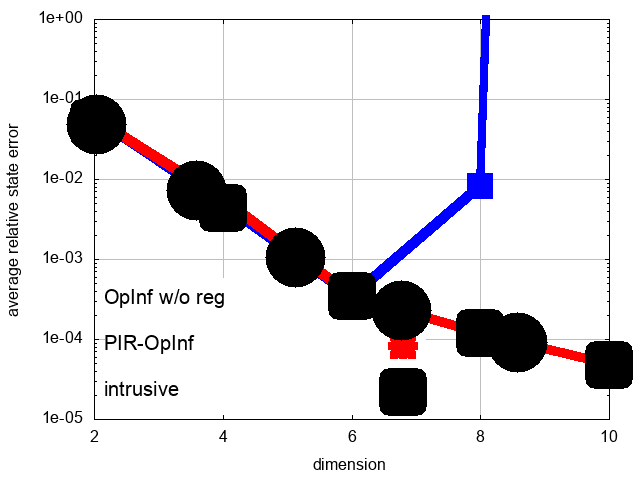}}}
  \label{fig: sample-test}}
    \hfill
    \subfloat[estimated stability radius]{{\Large\resizebox{0.48\columnwidth}{!}{\input{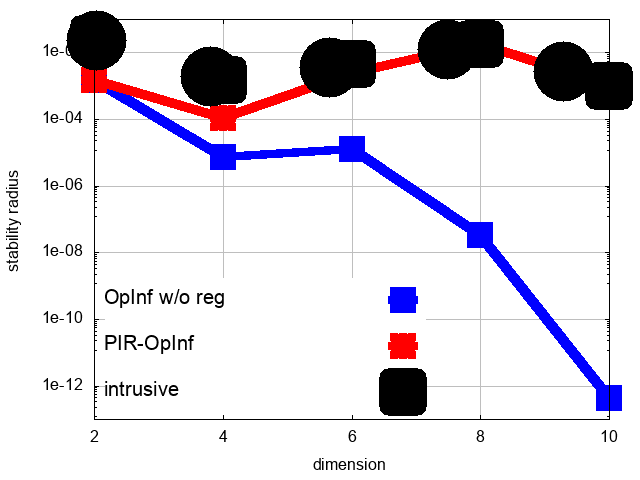}}}
    \label{fig: sample-rho}}
    \caption{Synthetic example: The model obtained with the proposed PIR-OpInf shows stable behavior, in contrast to OpInf without regularization, and achieves a comparable test error as intrusive model reduction. The estimated stability radius \eqref{eq:PIOpInf:Radius} of the PIR-OpInf model is orders of magnitude larger than the estimated stability radius of the OpInf model without regularization, which is in agreement with the aim of the proposed regularizer to penalize models with low stability radii.}
  \end{figure}

  Figure~\ref{fig: sample-test} shows the test error
  \begin{equation}
  e_{\text{test}} = \sum_{i = 1}^{\Mtest}\frac{\|\bfV \bar{\bfX}^{\text{test}}(\mutest_i) - \bfX^{\text{test}}(\mutest_i)\|_F}{\|\bfX^{\text{test}}(\mutest_i)\|_F}\,,
  \label{eq:NumExp:Syn:TestError}
  \end{equation}
  where $\bar{\bfX}^{\text{test}}(\mutest_i)$ is the trajectory obtained at test parameter $\mutest_i$ with the corresponding test input trajectory and test initial condition with either PIR-OpInf, OpInf without regularization, or intrusive model reduction. In contrast to OpInf without regularization, PIR-OpInf shows stable behavior and yields accurate predictions even for dimensions $\nr > 6$ in this example. Figure~\ref{fig: sample-rho} shows the stability radius $\hat{\rho}$ define in \eqref{eq:PIOpInf:Radius} for PIR-OpInf, OpInf without regularization, and intrusive model reduction. The stability radius of the model obtained with PIR-OpInf is larger than the stability radius of OpInf without regularization, which numerically demonstrates that the proposed physics-informed regularizer indeed induces a stability bias.

  \subsection{Burgers' equation}
  \label{sec:NumExp:Burgers}
  We consider the parameterized Burgers' equation
  \begin{equation*}
      \frac{\partial x}{\partial t}(\omega, t; \mu) = \mu \frac{\partial^2 x}{\partial^2 \omega}(\omega, t; \mu) - x(\omega, t; \mu) \frac{\partial x}{\partial \omega}(\omega, t; \mu)
  \end{equation*}
  with spatial coordinate $\omega \in (0, 1)$, time $t \in [0, 1]$, and viscosity $\mu \in [10, 100]$. Dirichlet boundary conditions $x(0, t; \mu) = u(t),\: x(1, t; \mu) = 0$ are imposed, with input $u: [0, 1] \to \mathbb{R}$. The equation is discretized in space with finite differences on an equidistant grid in $[0, 1]$ with $N = 128$ grid points. Time is discretized with the explicit Euler method with time-step size $\delta t = 10^{-4}$.

  \subsubsection{Problem setup}
  For each of the $M = 10$ training parameters $\mu = \{10, 20, 30, \dots, 100\}$, we derive a single input trajectory $\bfUbasis_1(\mu)$, with entries uniformly sampled in $[0, 2]$, and an initial condition $\bfx_1(\mu) = \boldsymbol{0}$. Thus, $\Mbasis = 1$. The corresponding state trajectories are $\bfXbasis_1(\mu_1), \dots, \bfXbasis_1(\mu_M)$.
  A basis matrix $\bfV \in \mathbb{R}^{\nh \times \nr}$ is then constructed from the corresponding snapshots as described in Section~\ref{section:irom}.
  Furthermore, we sample $\Mtrain = 10$ training inputs $\bfU_1(\mu), \dots, \bfU_{\Mtrain}(\mu)$ for each training parameter $\mu \in \{10, \dots, 100\}$. To generate the initial conditions $\bfx_{1,0}(\mu), \dots, \bfx_{\Mtrain,0}(\mu)$ for each training parameter $\mu \in \{10, \dots, 100\}$, we sample $\nr$-dimensional random vectors $\boldsymbol r_{1}(\mu), \dots, \boldsymbol r_{\Mtrain}(\mu)$ with independent entries uniformly distributed in $[0, 1]$ and set $\bfx_{i,0}(\mu) = \bfV\boldsymbol r_i(\mu)$ for $i = 1, \dots, \Mtrain$. We apply parameter selection as in Section~\ref{sec:NumExp:Syn} to find regularization parameters for each dimension $\nr \in \{2, \dots, 10\}$. Furthermore, we construct an operator-inference model obtained without regularization and a reduced model with intrusive model reduction.

  For comparison purposes, we also construct from the same training data an operator-inference model with the regularization proposed in \cite{mcquarrie2020data, doi:10.2514/1.J058943}, which is Tikhonov regularization that regularizers the Frobenius norms of linear, quadratic, and input operators together, rather than only the norm of the quadratic operator as in the proposed PIR-OpInf. We refer to this approach as T-OpInf in the following. The same regularization parameter-selection procedure is applied as for PIR-OpInf.

  \subsubsection{Results for PIR-OpInf}
  \label{sec:NumExp:Burgers:Results}
  We test the models at $\Mtest = 7$ test parameters that are equidistantly distributed in the parameter domain $\Dcal$. For each test parameter $\mu \in \{\mutest_1, \dots, \mutest_{\Mtest}\}$, we generate $\MtestIn = 5$ input trajectories $\bfU_1^{\text{test}}(\mu), \; \dots, \; \bfU_{\MtestIn}^{\text{test}}(\mu)$ and the corresponding test state trajectories $\bfX^{\text{test}}_1(\mu), \dots, \bfX^{\text{test}}_{\MtestIn}(\mu)$. The test initial conditions are the same as the training initial conditions, i.e, $\bfx^{\text{test}}_{i,0}(\mu) = \bfx_{i,0}(\mu)$, for $i = 1, \dots, \MtestIn$. The test error is then given by
  \begin{equation}
      e_{\text{test}} = \sum_{i = 1}^{\Mtest} \sum_{j = 1}^{\MtestIn} \frac{\|\bfV \bar{\bfX}_j^{\text{test}}(\mutest_i) - \bfX_j^{\text{test}}(\mutest_i)\|_F}{\|\bfX_j^{\text{test}}(\mutest_i)\|_F}\,,
      \label{sec:Comp:TestErr}
  \end{equation}
  where the trajectories $\bar{\bfX}_j^{\text{test}}(\mutest_i)$ for $i = 1, \dots, \Mtest$ and $j = 1, \dots, \MtestIn$ are obtained from either operator inference without regularization, PIR-OpInf, T-OpInf, or intrusive model reduction.

  \begin{figure}
    \centering
    \subfloat[input domain $(0,2)$ ]{{\Large\resizebox{0.45\columnwidth}{!}{\input{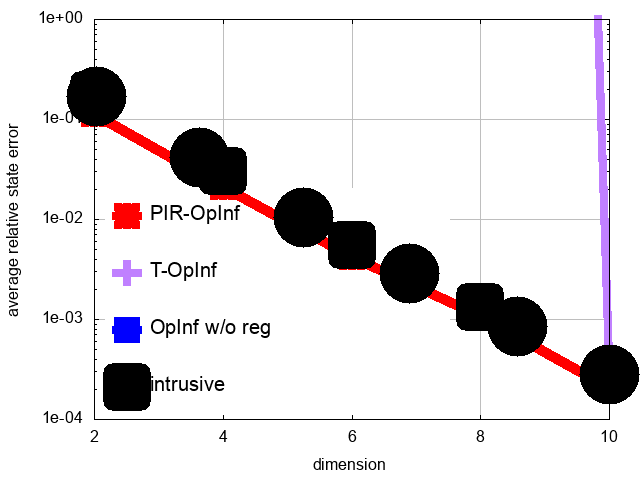}}}
    \label{fig:burgers-testU2}}
    \hfill
     \subfloat[input domain $(0,2)$, T-OpInf with post-processing]{{\Large\resizebox{0.45\columnwidth}{!}{\input{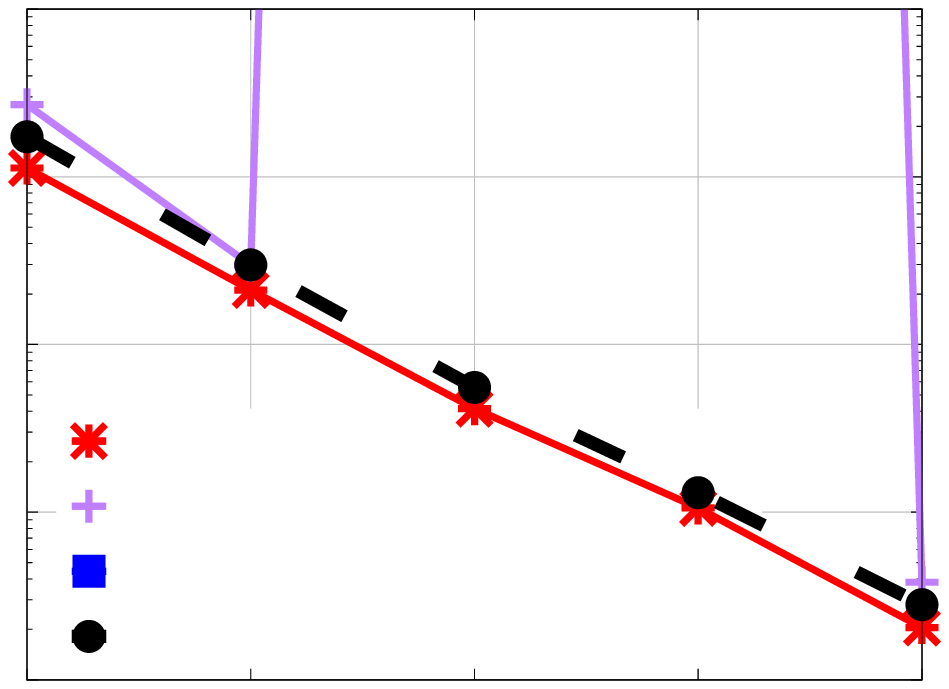}}}
    \label{fig:burgers-testU2post}}

    \subfloat[input domain $(0, 3)$]{{\Large\resizebox{0.45\columnwidth}{!}{\input{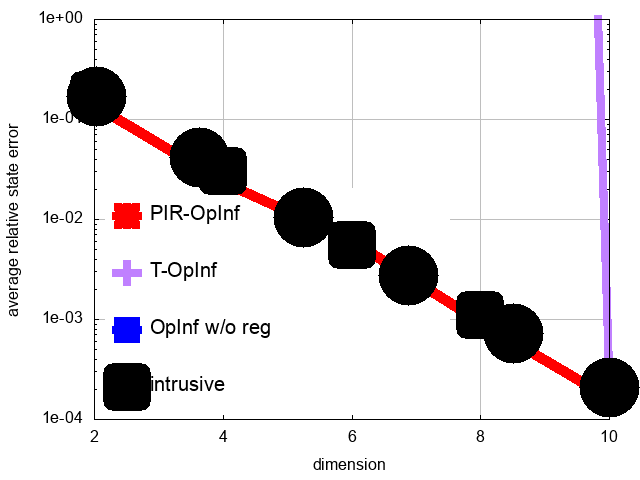}}}
    \label{fig:burgers-testU3}}
    \hfill
    \subfloat[input domain $(0, 3)$, T-OpInf with post-processing]{{\Large\resizebox{0.45\columnwidth}{!}{\input{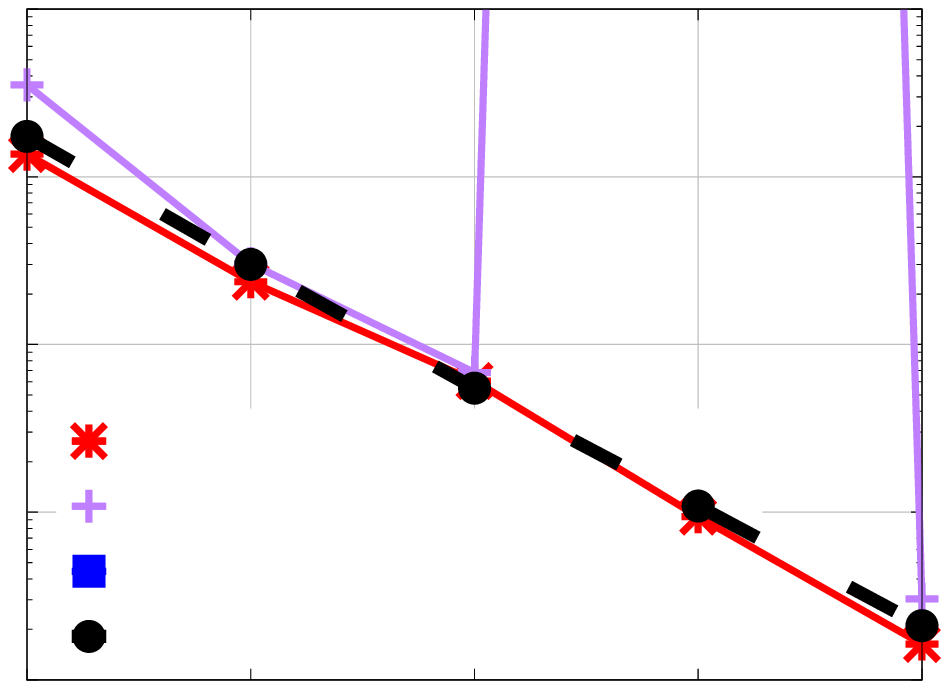}}}
    \label{fig:burgers-testU3post}}

    \subfloat[input domain $(0, 4)$]{{\Large\resizebox{0.45\columnwidth}{!}{\input{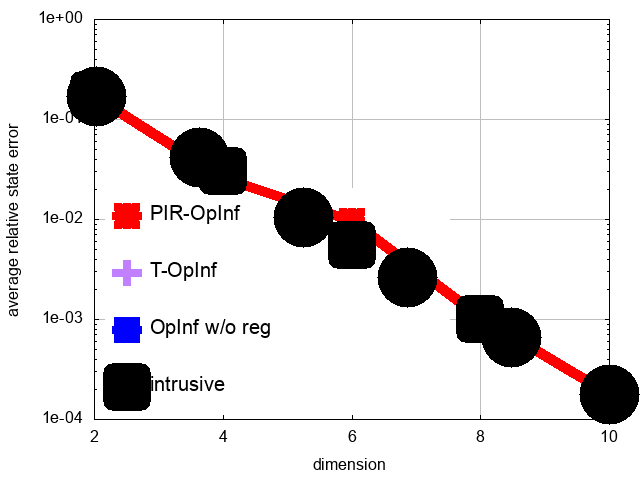}}}
    \label{fig:burgers-testU4}}
    \hfill
     \subfloat[input domain $(0, 4)$, T-OpInf with post-processing]{{\Large\resizebox{0.45\columnwidth}{!}{\input{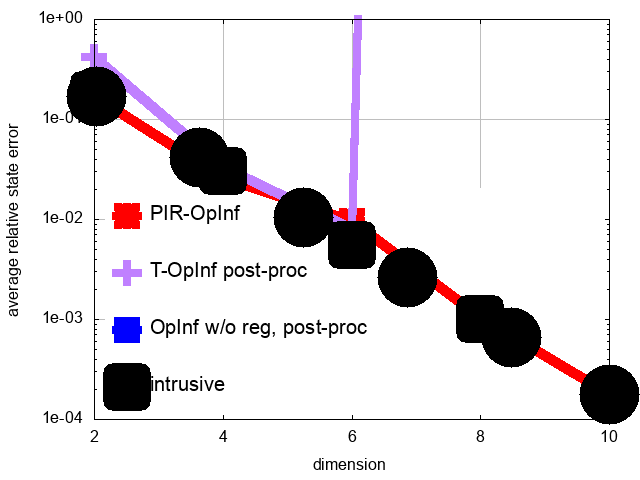}}}
    \label{fig:burgers-testU4post}}

    \caption{Burgers' equation: The proposed PIR-OpInf leads to models that are stable for a large range of inputs, which is in contrast to OpInf without regularization and OpInf with Tikhonov regularization (T-OpInf). The results also show that applying the same post-processing as for PIR-OpInf (cf.~Section~\ref{sec:PIOpInf:Reg}) to T-OpInf has little effect on the stability of the learned models, which indicates that indeed the proposed regularizer in PIR-OpInf is responsible for obtaining stabler models.}
    \label{fig:burgers-testU2-4all}
  \end{figure}

  Figure~\ref{fig:burgers-testU2}, Figure~\ref{fig:burgers-testU3}, and Figure~\ref{fig:burgers-testU4} show the test error \eqref{sec:Comp:TestErr} for test inputs with entries sampled uniform from the domains $[0, 2]$, $[0, 3]$, and $[0, 4]$, respectively. In all cases, PIR-OpInf shows stable behavior, whereas OpInf without regularization leads to numerical instabilities. Even T-OpInf with Tikhonov regularization shows unstable behavior for many dimensions $\nr$. To separate the effect of the regularization from the effect of the post-processing (cf.~Section~\ref{sec:PIOpInf:Reg}), we apply the same post-processing as for PIR-OpInf to T-OpInf. The corresponding results in Figure~\ref{fig:burgers-testU2post}, Figure~\ref{fig:burgers-testU3post}, and Figure~\ref{fig:burgers-testU4post} show that post-processing helps to stabilize T-OpInf as well; however, as the range of the inputs increases, a similarly unstable behavior as in the case without post-processing is obtained. Thus, the results indicate that penalizing the quadratic term via the proposed regularizer is responsible for achieving stabler models, rather than the post processing or penalizing both the linear and the quadratic term together as in Tikhonov regularization, which is in agreement with the theoretical motivation outlined in Section~\ref{sec:PIOpInf:Stability}.

\begin{figure}
    \centering
    \subfloat[T-OpInf without post-processing]{{\Large\resizebox{0.45\columnwidth}{!}{\input{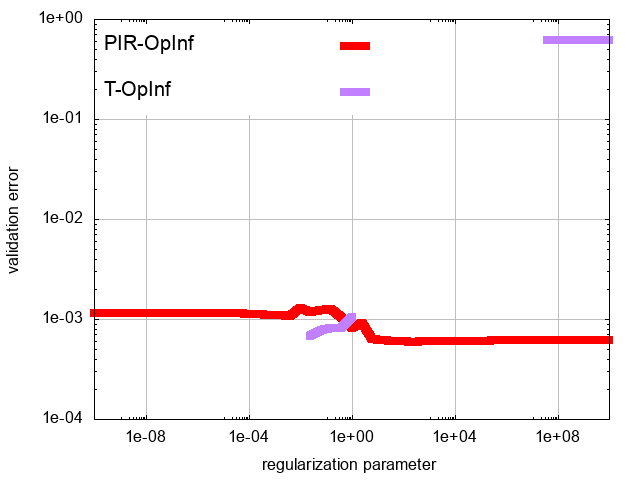}}}
    \label{fig:burgersdim8-cv}}
    \hfill
    \subfloat[T-OpInf with post-processing]{{\Large\resizebox{0.45\columnwidth}{!}{\input{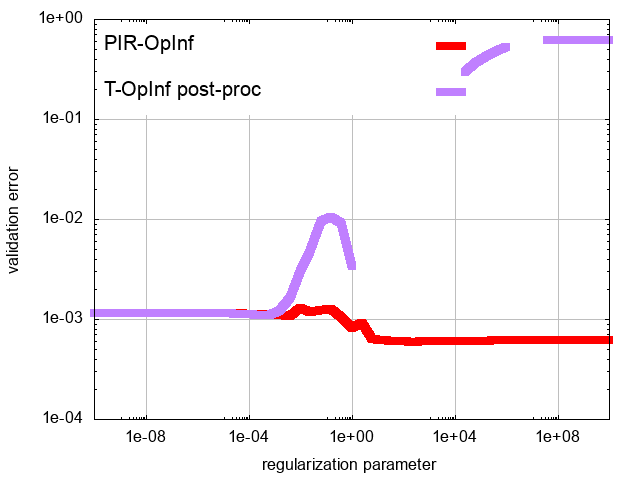}}}
    \label{fig:burgersdim8-cv-post}}

    \caption{Burgers' equation: The validation error (objective of \eqref{sec:Comp:ParamSelOpti}) of T-OpInf grows quickly with the regularization parameter, which means that the regularization has a negative effect on the model accuracy. In contrast, the validation error corresponding to the proposed PIR-OpInf is less sensitive to the regularization parameter, which means that large regularization parameters can be chosen---imposing a stronger stability bias---without deteriorating the model accuracy.}
    \label{fig:burgersdim8-cvall}
  \end{figure}

  Consider now Figure~\ref{fig:burgersdim8-cvall} that shows the validation error, i.e., the objective of \eqref{sec:Comp:ParamSelOpti}, of the parameter-selection procedure versus the regularization parameter for dimension $\nr = 8$ for PIR-OpInf and T-OpInf. Independent of whether post-processing is applied to T-OpInf (Figure~\ref{fig:burgersdim8-cv-post}) or not (Figure~\ref{fig:burgersdim8-cv}), the error of T-OpInf grows quickly as the regularization parameter is increased. Thus, if a small regularization parameter is chosen, the validation error of T-OpInf is small but it also means that no regularization is induced. If instead the regularization parameter is large, then there might be a stability bias but at the same time it leads to a distinct increase of the model error. In contrast, the curves corresponding to PIR-OpInf show that the validation error is small for moderately sized regularization parameters, where a stability bias is induced without leading to a deterioration of the model accuracy.

  \begin{figure}
    \centering
    \subfloat[estimated stability radius]{{\Large\resizebox{0.45\columnwidth}{!}{\input{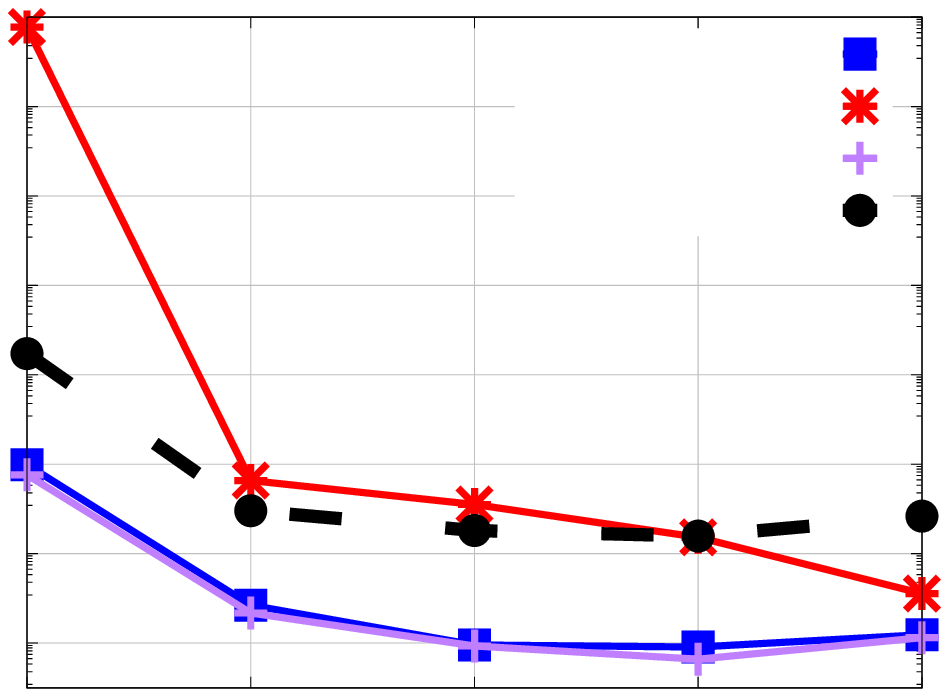}}}
    \label{fig:burgers-rho-nee}}
    \subfloat[validation error for $\nr = 2$]{{\Large\resizebox{0.45\columnwidth}{!}{\input{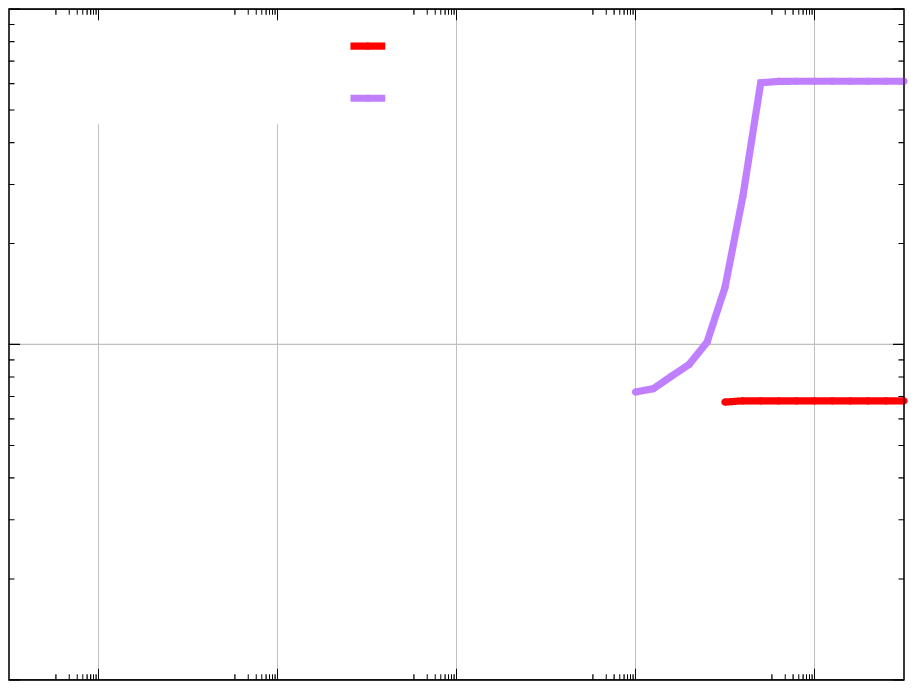}}}
    \label{fig:burgersdim2-cv}}
    \caption{Burgers' equation: The models learned with the proposed PIR-OpInf have larger estimated stability radii than models obtained without regularization and with Tikhonov regularization (T-OpInf). The stability radius for $\nr = 2$ is high for PIR-OpInf models because only large regularization parameters $\lambda^* = 10^6$ lead to stable behavior (see (b)), which forces the norm of the quadratic term to be close to zero and thus increases the stability radius by a large amount. Notice that the PIR-OpInf model at $\nr = 2$ achieves a similar accuracy as intrusive model reduction; cf.~Figure~\ref{fig:burgers-testU2-4all}.}
    \label{fig:burgers-rho}
  \end{figure}

  Figure~\ref{fig:burgers-rho-nee} shows the estimated stability radius \eqref{eq:PIOpInf:Radius} for various models. The results indicate that PIR-OpInf achieves a larger stability radius than the models obtained with T-OpInf and OpInf without regularization. At dimension $\nr = 2$, the stability radius of PIR-OpInf is large because only large regularization parameters lead to stable behavior; see Figure~\ref{fig:burgersdim2-cv}. The regularization parameter is chosen so large that the quadratic term is close to 0, which explains the high estimated stability radius.

  \subsubsection{Results for SPIR-OpInf}
\begin{figure}
    \centering
    \subfloat[input domain $(0,2)$ ]{{\huge\resizebox{0.32\columnwidth}{!}{\input{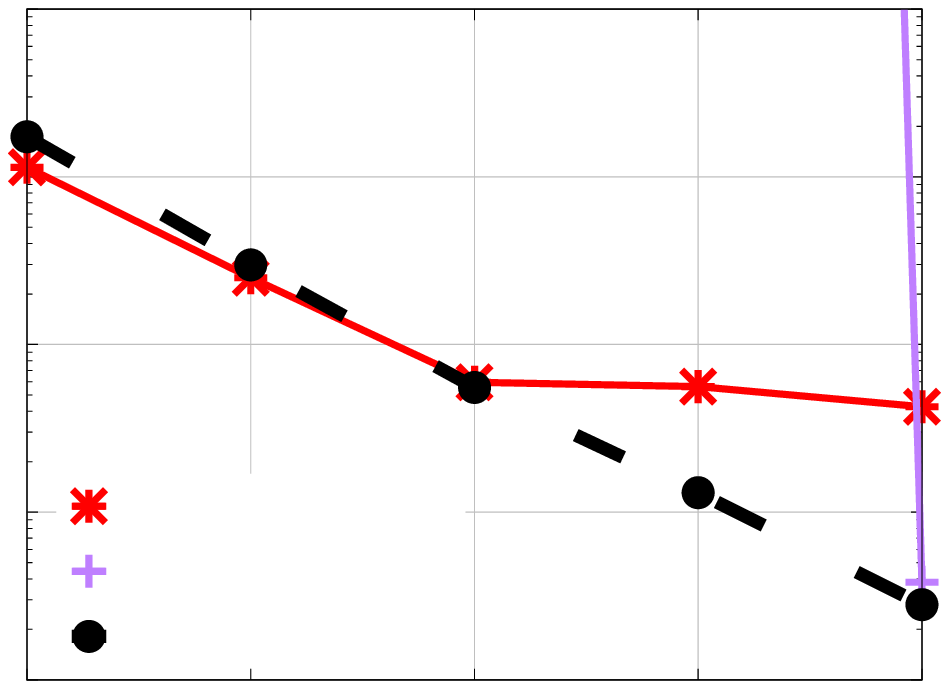}}}
    \label{fig:burgers-testU2spd}}
    \hfill
    \subfloat[input domain $(0, 3)$]{{\huge\resizebox{0.32\columnwidth}{!}{\input{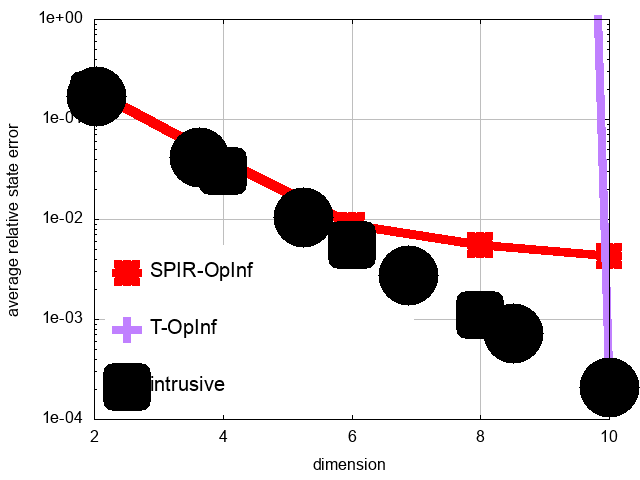}}}
    \label{fig:burgers-testU3spd}}
    \hfill
    \subfloat[input domain $(0, 4)$]{{\huge\resizebox{0.32\columnwidth}{!}{\input{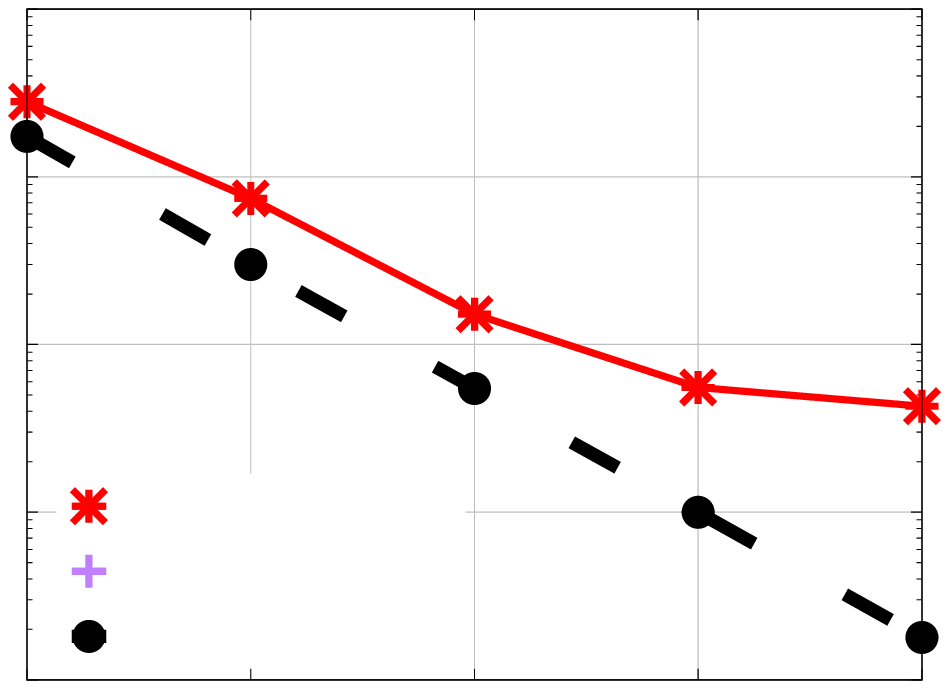}}}
    \label{fig:burgers-testU4spd}}
    \caption{Burgers' equation: Imposing symmetry and definiteness onto the linear operator with the proposed SPIR-OpInf \eqref{eq:PIOpInf:ObjSPIR} leads to stable models in this experiment; however, the additional constraints lead to a lower accuracy than PIR-OpInf that include the proposed regularizer but no constraints on the linear operator (cf.~Figure~\ref{fig:burgers-testU2-4all}).}
    \label{fig:burgers-testU2-4spdall}
  \end{figure}

\begin{figure}
    \centering
    \subfloat[estimated stability radius]{{\Large\resizebox{0.45\columnwidth}{!}{\input{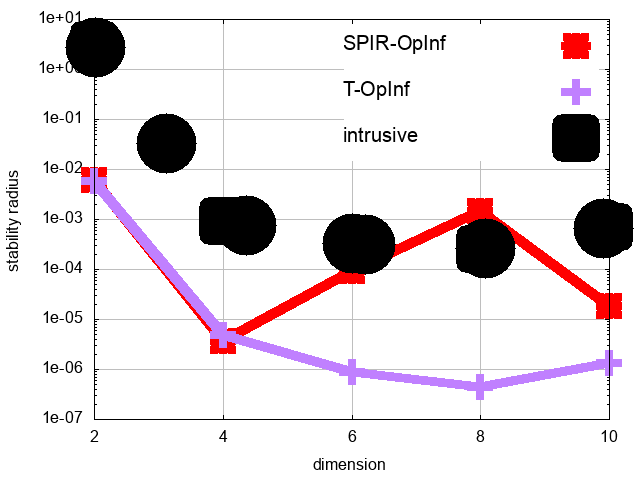}}}
    \label{fig:burgers-rho-see}}
    \subfloat[validation error for dimension $\nr = 2$]{{\Large\resizebox{0.45\columnwidth}{!}{\input{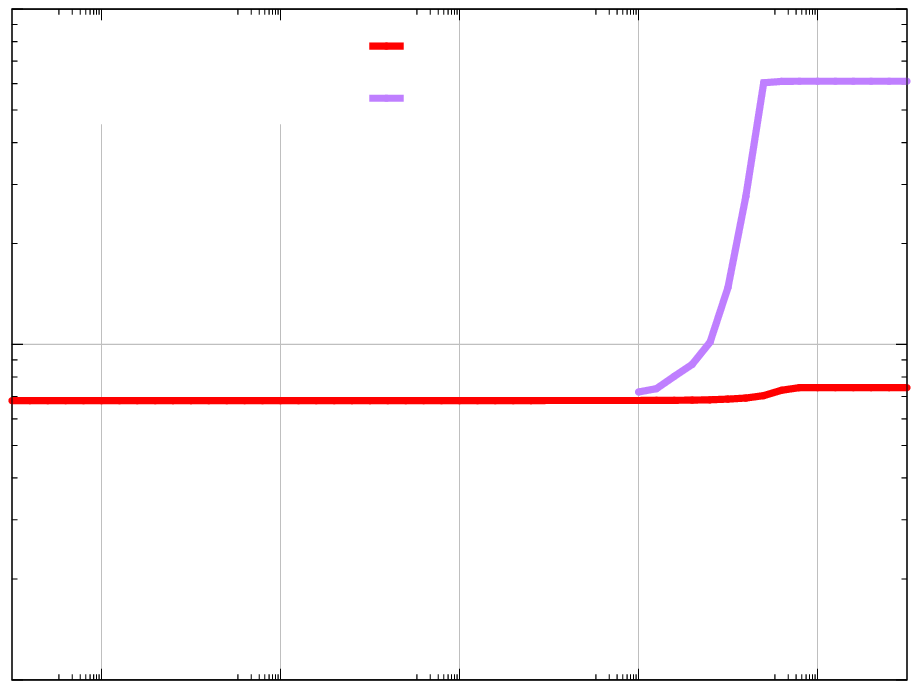}}}
    \label{fig: burgersdim2-cv-see}}

    \caption{Burgers' equation: Models learned with SPIR-OpInf have a larger estimated stability radius than models learned with Tikhonov regularization for dimensions $\nr > 4$ in this example. The stability radius of the SPIR-OpInf model is similar to the stability radius of the T-OpInf model for $\nr = 2$ because a small regularization parameter is chosen as shown in plot (b).}
    \label{fig:burgers-rho-seeall}
  \end{figure}

  Figure~\ref{fig:burgers-testU2-4spdall} shows the test error \eqref{sec:Comp:TestErr} corresponding to SPIR-OpInf, which imposes symmetry and definiteness onto the linear operator. Stable behavior is obtained in all cases; however, a leveling off of the error as the dimension increases indicates that restricting to symmetric operators in this example is limiting the accuracy. The estimated stability radii of the SPIR-OpInf models are compared to the stability radii of T-OpInf models in Figure~\ref{fig:burgers-rho-see}. For large dimensions $\nr > 4$, the estimated stability radii of the SPIR-OpInf models is larger than the stability radii of the T-OpInf models. For small dimensions $\nr \leq 4$, the stability radii of SPIR-OpInf and T-OpInf models is large. This is reflected by the small regularization parameter selected by the proposed parameter selection procedure, which selects $\lambda \approx 3 \times 10^{-7}$ for $\nr = 2$.

  \subsection{Reactive-diffusion problem}
  \label{sec:NumExp:React}
  Consider the parameterized reaction-diffusion equation
  \begin{equation}
  \frac{\partial}{\partial t}x(\bfXi, t; \mu) =
     \Delta x(\bfXi, t; \mu) + s(\bfXi)u(t) + g(x(\bfXi, t; \mu))\,,
  \end{equation}
  with the spatial coordinate $\bfXi = [\xi_1 \; \xi_2]^T \in [0, 1]^2$. We impose homogeneous Neumann boundary conditions. The parameter domain is $\Dcal = [1, 1.5]$ and end time is $T = 20$. The source is $s(\bfXi) = 10^{-1} \operatorname{sin}(2 \pi \xi_1) \operatorname{sin}(2 \pi \xi_2)$. The non-linear term is
  \begin{equation*}
      g(x(\bfXi, t; \mu)) = -(a \operatorname{sin}(\mu) + 2) \operatorname{exp}(-\mu^2 b) \bigg(1 + (\mu c) x  + \frac{{(\mu c)}^2}{2!} x^2 \bigg)
  \end{equation*}
  which is the second-order Taylor approximation of the source term used in \cite{peherstorfer2019sampling}, with $a = 0.1, b = 2.7$ and $c = 1.8$. We discretize in space with a mesh width of $h = 1/12$ and finite difference and in time with explicit Euler with time-step size $\delta t = 10^{-2}$. The dimension of the high-dimensional model is $\nh = 144$.

  The training parameter set contains the $M = 10$ equidistant points in the parameter domain $\Dcal$. To construct the reduced space, we take a single $\Mbasis = 1$ input trajectory for each training parameter, where the inputs are sampled uniformly in $[0, 1]$. The initial condition is zero. The corresponding trajectories $\bfXbasis_1(\mu_1), \dots, \bfXbasis_1(\mu_M)$ are used to construct a POD basis.
  For each of the $M$ training parameters, we sample $\Mtrain = 10$ input trajectories with entries uniformly in $[0, 1]$. The regularization parameters are selected via our selection procedure described in Section~\ref{sec:Comp:CV} by sweeping over the 51 logarithmically equidistant points in $[10^{-10}, 10^{10}]$, which is the same setup as in the previous experiments. We test models for $\Mtest = 7$ test parameters that are equidistantly distributed in the parameter domain $\Dcal$. For every test parameter $\mutest_1, \dots, \mutest_{\Mtest}$, we generate a single input trajectory $\bfU^{\text{test}}$, whose entries randomly selected via a uniform distribution in $[0, 1]$, and initial condition $\bfx_0^{\text{test}} = \boldsymbol{0}$. The test error is then
  \begin{equation}
      e_{\text{test}} = \sum_{i = 1}^{\Mtest} \frac{\|\bfV \bar{\bfX}^{\text{test}}(\mutest_i) - \bfX^{\text{test}}(\mutest_i)\|_F}{\|\bfX^{\text{test}}(\mutest_i)\|_F}\,,
      \label{sec:Comp:TestErrReact}
  \end{equation}
  where $\bar{\bfX}^{\text{test}}(\mutest_i)$ is the predicted trajectory at parameter $\mutest_i$ by either the PIR-OpInf, T-OpInf, OpInf without regularization, or intrusive model reduction.

\begin{figure}
    \centering
    \subfloat[test error \eqref{sec:Comp:TestErrReact}]{{\Large\resizebox{0.45\columnwidth}{!}{\input{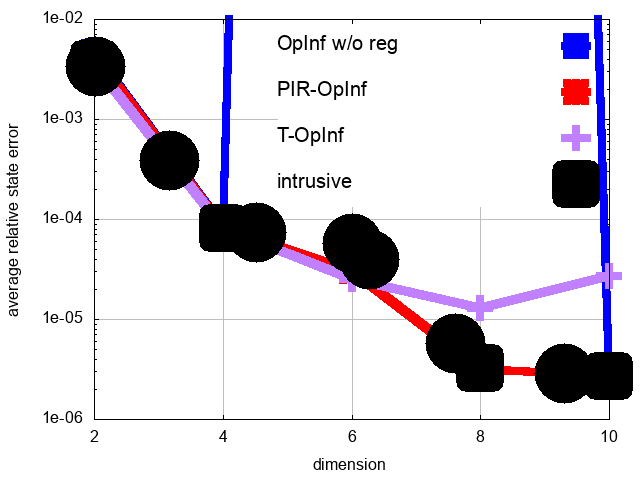}}}
  \label{fig: react-nee}}
    \hfill
    \subfloat[estimated stability radius \eqref{eq:PIOpInf:Radius}]{{\Large\resizebox{0.45\columnwidth}{!}{\input{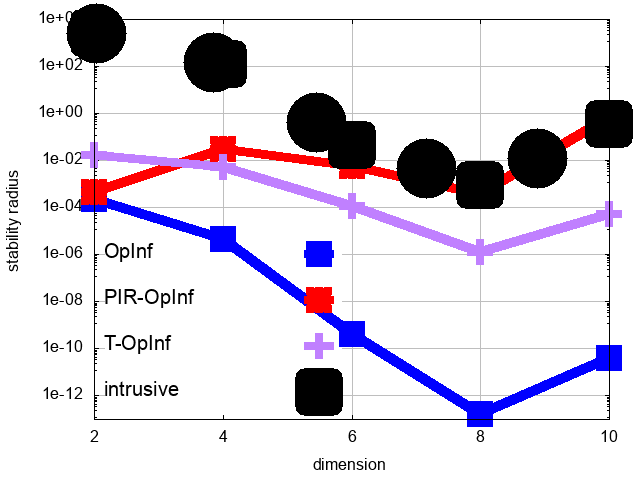}}}
  \label{fig: react-rho-nee}}
    \caption{Reaction-diffusion problem: The PIR-OpInf model shows stable behavior in this experiment. In contrast to Tikhonov regularization (T-OpInf), the PIR-OpInf model achieves an accuracy close to intrusive model reduction even for larger $\nr > 6$ dimensions. The estimated stability radius of the PIR-OpInf model is orders of magnitude higher than the stability radius of the T-OpInf model in this experiment.}
    \label{fig:ReactDiffAll}
  \end{figure}

\begin{figure}
    \centering
    \subfloat[test error \eqref{sec:Comp:TestErrReact}]{{\Large\resizebox{0.45\columnwidth}{!}{\input{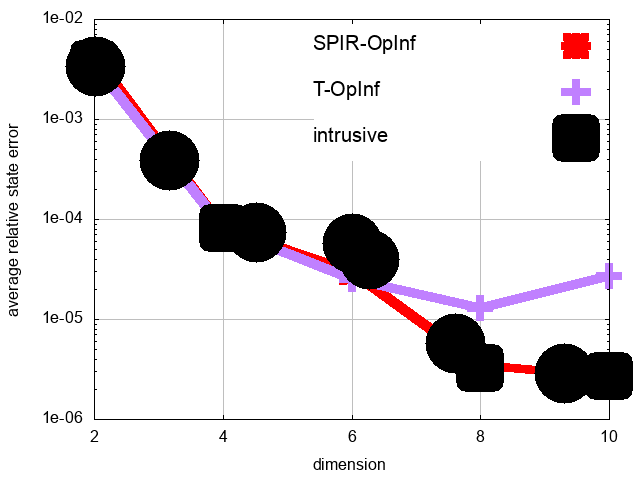}}}
     \label{fig: react-see}}
    \hfill
    \subfloat[estimated stability radius \eqref{eq:PIOpInf:Radius}]{{\Large\resizebox{0.45\columnwidth}{!}{\input{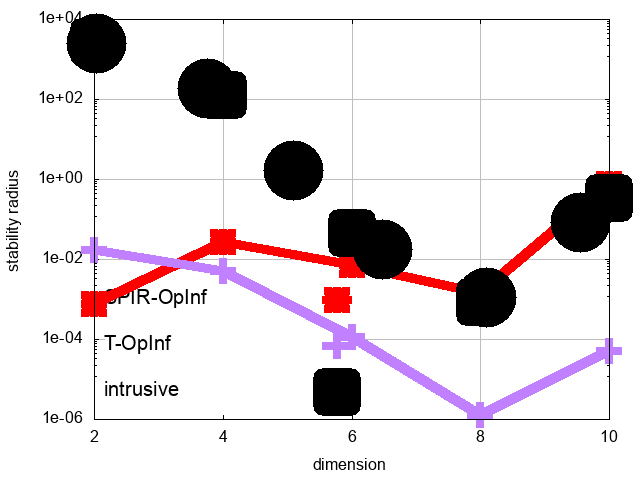}}}
    \label{fig: react-rho-see}}

    \caption{Reaction-diffusion problem: Constraining the linear inferred operator to be symmetric and definite with SPIR-OpInf leads to models with comparable accuracy and stability radius as PIR-OpInf in this example.}
    \label{fig:ReactDiffSPDAll}
  \end{figure}

  Figure~\ref{fig: react-nee} shows the test error \eqref{sec:Comp:TestErrReact}. The results indicate that OpInf without any regularization becomes unstable quickly. In contrast, T-OpInf and PIR-OpInf provide stable approximations. However, whereas the Tikhonov regularization in T-OpInf leads to a loss of accuracy at higher dimensions, the proposed physics-informed regularizer used by PIR-OpInf achieves errors that are comparable to intrusive model reduction. This is in agreement with the estimated stability radius shown in Figure~\ref{fig: react-rho-nee}, where PIR-OpInf achieves an orders of magnitude larger stability radius at higher dimensions $\nr$ than T-OpInf and OpInf without regularization. Similar results are obtained with SPIR-OpInf, where symmetry and definiteness are imposed, as shown in Figure~\ref{fig:ReactDiffSPDAll}.

  \section{Conclusions}
  \label{sec:Conc}
  Learning models from data is an ever more important task in science and engineering. It is increasingly recognized that physical insights need to be incorporate together with data to learn truly predictive models \cite{doi:10.1098/rsta.2016.0153,Willcox2021}. In this spirit, we proposed a regularizer that explicitly leverages the quadratic model form, which in turn is imposed by the underlying physics, to penalize unstable models learned with operator inference. We also showed that additional physical insights in the form of structure of the linear dynamics can be imposed on the operator-inference models via constraints. In our experiments, operator inference with the proposed physics-informed regularizer and structure preservation outperforms operator inference without regularization and operator inference with Tikhonov regularization in terms of stability and accuracy. Thus, our results provide evidence of the importance of combining physical insights and data for deriving predictive models in science and engineering.

  \section*{Acknowledgements}
  The first and third author acknowledge partially supported by US Department of Energy, Office of Advanced Scientific Computing Research, Applied Mathematics Program (Program Manager Dr. Steven Lee), DOE Award DESC0019334,
  and by the National Science Foundation under Grant No.~1901091 and under Grant No.~1761068.

  \bibliography{References.bib}

\bibliographystyle{abbrv}

\end{document}